\documentclass[11pt,letterpaper,reqno]{amsart}
\usepackage[T1]{fontenc}
\usepackage[utf8]{inputenc}
\usepackage{lmodern,microtype}
\usepackage{amsmath,amssymb,mathtools,bm,mathrsfs}
\usepackage[margin=1.25in,headheight=14pt,headsep=6mm]{geometry}
\usepackage{enumitem}
\usepackage[hidelinks]{hyperref}
\hypersetup{
	pdftitle={Invariant measures and uniform small-forcing limits for radially stabilized Hamiltonian transport equations},
	pdfauthor={Yong Chen and Hongjun Gao},
	pdfsubject={Statistical stability and invariant measures for stochastic transport dynamical systems}
}
\allowdisplaybreaks[2]
\numberwithin{equation}{section}
\newtheorem{theorem}{Theorem}[section]
\newtheorem{lemma}[theorem]{Lemma}
\newtheorem{proposition}[theorem]{Proposition}
\newtheorem{corollary}[theorem]{Corollary}
\theoremstyle{definition}
\newtheorem{assumption}[theorem]{Assumption}
\newtheorem{definition}[theorem]{Definition}
\theoremstyle{remark}
\newtheorem{remark}[theorem]{Remark}
\newtheorem{example}[theorem]{Example}
\newcommand{\R}{\mathbb{R}}
\newcommand{\T}{\mathbb{T}}

\newcommand{\E}{\mathbb{E}}
\newcommand{\PP}{\mathbb{P}}
\newcommand{\dd}{\,\mathrm{d}}
\newcommand{\HS}{\mathcal{L}_2}
\newcommand{\ip}[2]{\langle #1,#2\rangle}
\newcommand{\norm}[1]{\lVert #1\rVert}
\newcommand{\PhiV}{\Phi_p}
\newcommand{\cI}{\mathcal{I}}
\newcommand{\cJ}{\mathcal{J}}

\DeclareMathOperator{\Law}{Law}
\DeclareMathOperator{\Lip}{Lip}
\DeclareMathOperator{\Tr}{Tr}

\title[Stationary scaling for Hamiltonian transport]{Invariant measures and uniform small-forcing limits for radially stabilized Hamiltonian transport equations}
\author[Yong Chen]{Yong Chen}
\thanks{\raggedright Yong Chen: Department of Mathematics, Zhejiang Sci-Tech
	University, Hangzhou 310018, China. E-mail: \texttt{youngchen329@126.com}.}
\author[Hongjun Gao]{Hongjun Gao}
\thanks{\raggedright Hongjun Gao: School of Mathematics, Southeast University,
	Nanjing 211189, China. E-mail: \texttt{hjgao@seu.edu.cn}.}
\date{}
\subjclass[2020]{Primary 60H15; Secondary 35Q35, 37L40, 60G10}
\keywords{Hamiltonian transport equations, radial It\^o noise, invariant measures, small-forcing limit, Cauchy law, Euler--Poincar\'e equations}
\begin{document}
	
\begin{abstract}
We study invariant measures and small-forcing limits for inviscid stochastic transport equations whose deterministic drifts arise from Hamiltonian, Lie--Poisson, or Euler--Poincar\'e dynamics. A common scalar It\^o noise acts radially on the state, while an independent additive noise has amplitude $\eta$. Since the quadratic transport drift may lose one spatial derivative, we develop a two-space well-posedness theory in which trajectories are controlled in a strong Hilbert space and stability is measured in a weaker space. A structure-preserving regularization recovers the high-order energy identity, and a concave fractional energy yields global bounds uniform in $\eta$. For every positive forcing amplitude, invariant measures exist and every such measure satisfies the same stationary scaling estimate. When the additive covariance is nonzero, the entire family of rescaled invariant measures converges, with an explicit rate, to a Gaussian-ratio Cauchy law. We also prove convergence of fractional moments, fractional Wasserstein distances, and stationary path laws. The results apply to active scalars, including two-dimensional Euler and surface quasi-geostrophic equations, ideal magnetohydrodynamics, smoothing-inertia EPDiff equations, and the Camassa-Holm, Degasperis-Procesi, and Hunter-Saxton equations.
\end{abstract}
	\maketitle

\section{Introduction}\label{sec:intro}

\subsection{Background and motivation}
Inviscid transport equations often admit a geometric formulation that explains
both their conservation laws and their leading-order cancellations. The active
scalar models considered here have Lie--Poisson formulations. Ideal
magnetohydrodynamics(MHD) is associated with a semidirect-product
Lie--Poisson bracket. The Euler--Poincar\'e equations on the diffeomorphism
group (EPDiff), together with the Camassa--Holm and Hunter--Saxton equations,
arise from Euler--Poincar\'e (EP) reduction. Camassa--Holm and
Degasperis--Procesi also have bi-Hamiltonian formulations. We refer to these
related classes collectively as \emph{Hamiltonian transport equations}; see
\cite{EM,AK,HMR,CH,KM03,DHH}.

The geometric structure identifies distinguished energies and cancellations,
but it does not by itself provide global control of high Sobolev norms. This is
especially important for transport equations, whose nonlinear drift typically
loses one spatial derivative. Stochastic perturbations create a second issue:
one must distinguish noise that changes the spatial characteristics from noise
that acts only on the amplitude of the state. The present paper concerns the
second type. A single scalar It\^o noise acts radially on the full state, while
an independent additive noise provides persistent forcing.

For dissipative stochastic partial differential equations, invariant measures
are commonly obtained from coercive energy inequalities and compactness.
Stochastic Navier--Stokes equations provide a standard example
\cite{FG,DPZ}. For inviscid Euler dynamics, statistical solutions have also
been constructed through fluctuation--dissipation or vanishing-viscosity limits
\cite{KS03,Kuksin04,KSbook,GHSV,BFerrari}. In those approaches, a dissipative
approximation supplies the compactness needed before the inviscid limit is
taken.

The long-time theory for the models considered here is less uniform. For the
surface quasi-geostrophic (SQG) equation, invariant-measure problems in
dissipative and inviscid regimes are studied in \cite{RZZ15,FS21}. Stationary
solutions for a modified SQG equation with transport noise appear in
\cite{LZ21}. For MHD, invariant measures and stationary limits have mainly
been obtained by adding viscosity, taming, or singular perturbations
\cite{FFGR17,Motyl22,Schenke22}. Stochastic Camassa--Holm equations have been
studied under several damping and noise mechanisms; see, for example,
\cite{APT}. For EPDiff, Degasperis--Procesi, and Hunter--Saxton, most of
the literature concerns geometric formulations, well-posedness, wave breaking,
and coherent structures
\cite{BauerEtAl20,HT16,TangEP23,CGDP15,ACCDP23,CCDP24,CGLDP25,HKP21}.

Our aim is to give a common stationary theory for these inviscid models under
radial It\^o stabilization. The same framework treats local well-posedness,
non-explosion, invariant measures, and the small-forcing limit. It also
separates the model-dependent part of the argument from the universal limiting
statistics: the geometric structure controls the nonlinear flow, whereas the
rescaled stationary limit depends only on the radial coefficient at the origin
and on the covariance of the additive noise.

\subsection{The stochastic model and its geometric input}
Let $V\Subset X$ be real separable Hilbert spaces. We consider
\begin{equation}\label{eq:original}
	\dd U_\eta+F(U_\eta)\dd t
	=\psi(\rho(U_\eta))U_\eta\dd W_t+\eta G\dd B_t,
	\qquad 0\le\eta\le1.
\end{equation}
Here $W$ is a real Brownian motion, $B$ is an independent cylindrical Brownian
motion on a separable Hilbert space $\mathfrak U$, and
$G\in\mathcal L_2(\mathfrak U,V)$. We write $G_V^*:V\to\mathfrak U$ for the
adjoint of $G$ and set
\[
Q_V=GG_V^*,\qquad
q_V=\|G\|_{\mathcal L_2(\mathfrak U,V)}^2.
\]
The parameter $\eta$ is the amplitude of the additive forcing. The functional
$\rho$ is one-homogeneous and measures a lower-order amplitude, while the drift
$F$ is quadratic. Depending on the equation, $U$ may be a vorticity, a
momentum, a velocity, or a product of transported fields.

The two Hilbert spaces have different roles. The strong space $V$ carries the
solution paths and the high-order energy. The weaker space $X$ is used for
stability, uniqueness, and weak convergence of probability measures. The drift
may satisfy only
\[
F:V\longrightarrow X,
\]
so the $V$-energy pairing is not automatically meaningful. For smooth states,
the formal deterministic contribution to the squared $V$-norm is
\begin{equation}\label{eq:intro-bv}
	b_V(u):=-2\langle F(u),u\rangle_V.
\end{equation}
Since $F(u)$ may be only $X$-valued, \eqref{eq:intro-bv} is first interpreted
on a smooth core and then recovered from regularized equations. The abstract
assumptions require this formal energy production to extend to a continuous
cubic functional on $V$ such that
\begin{equation}\label{eq:intro-bv-bound}
	|b_V(u)|\le C_V\rho(u)\|u\|_V^2.
\end{equation}
This estimate is the analytic form of the top-order cancellation supplied by
the geometric structure. The particular expression for $b_V$ is verified
separately for each model in Section~\ref{sec:models}.

The multiplicative noise in \eqref{eq:original} is radial: the same scalar
Brownian motion multiplies the entire state by
$\psi(\rho(U_\eta))$. This differs from Stratonovich transport noise and
stochastic Lie transport, which perturb the spatial characteristics
\cite{Holm15,FP}. Radial noise does not improve spatial regularity. Its role is
to modify the energy drift through the It\^o correction. The coefficient
$\psi$ is chosen so that this correction dominates the cubic term
\eqref{eq:intro-bv-bound} at large amplitude. We regard this feedback as an idealized amplitude-dependent
stabilization mechanism: it acts on the size of the state rather than on its
spatial derivatives. It is therefore analytically distinct from viscosity and transport noise.
The Hamiltonian formulation supplies cancellations for the deterministic
PDE; it does not imply that this radial It\^o perturbation is generated
by a stochastic Hamiltonian variational principle. The conclusions concern
this explicitly prescribed feedback forcing.

\subsection{Analytical strategy and stationary scaling}
The proof has four main steps. First, we construct local solutions by
structure-preserving regularization. The regularized equations satisfy uniform
high-order bounds, while their differences are controlled in $X$. Compactness
and stability then identify a limit on the original stochastic basis. Passing
to the limit in the regularized energy pairings recovers the $V$-energy
identity. This two-space argument is related to the regularization and
lower-topology methods used for stochastic transport and fluid equations
\cite{FGP,Flandoli11,GHV,BFM,CFH,CDG23,TY23,ART}.

Second, the recovered energy identity reveals the stabilization mechanism. If
$Y=\|U_\eta\|_V^2$ and $a(u)=\psi(\rho(u))$, then formally
\begin{equation}\label{eq:intro-energy}
	\dd Y=[b_V(U_\eta)+a(U_\eta)^2Y+\eta^2q_V]\dd t
	+2a(U_\eta)Y\dd W_t
	+2\eta\langle G_V^*U_\eta,\dd B_t\rangle_{\mathfrak U}.
\end{equation}
Apply It\^o's formula to $(1+Y)^p$ with $0<p<1/2$. The two contributions from
the radial noise combine into
\[
p\,a(u)^2Y(1+Y)^{p-2}\bigl[1-(1-2p)Y\bigr].
\]
This term is negative when $Y$ is large. If
\[
\psi(y)^2\ge b_0+b_1y,
\qquad (1-2p)b_1>C_V,
\]
then the negative quadratic-variation contribution dominates the transport
bound \eqref{eq:intro-bv-bound}. We obtain a coercive fractional Lyapunov
estimate and hence non-explosion. This mechanism is related to the
noise-induced stabilization results in
\cite{GHV,TY23,AKM,BMX,BagnaraS,LCcontrol}.

Third, the same Lyapunov estimate yields stationary compactness. Continuity of
the transition probabilities is proved on bounded subsets of $V$, with the
$X$ topology. This restricted continuity is enough to pass to the limit in
Krylov--Bogoliubov averages. A stationary truncation argument then extends the
Lyapunov bound from the constructed invariant measure to every invariant
probability. In the original variables, the estimate identifies $\eta$ as the
natural stationary amplitude.

Finally, the homogeneity of $F$ and $\rho$ determines the small-forcing
limit. For $\eta>0$, set $v_\eta=U_\eta/\eta$. Then
\begin{equation}\label{eq:intro-scaled}
	\dd v_\eta+\eta F(v_\eta)\dd t
	=\psi(\eta\rho(v_\eta))v_\eta\dd W_t+G\dd B_t.
\end{equation}
As $\eta\downarrow0$, the nonlinear drift vanishes and the radial coefficient
converges to $a_0=\psi(0)$. The limiting equation is
\[
\dd V=a_0V\dd W_t+G\dd B_t.
\]
Its stationary solution is a Gaussian variance mixture. More precisely,
\[
\nu_*=\Law\!\left(\frac{Z}{a_0|\zeta|}\right),
\qquad Z\sim\mathcal N(0,Q_V),
\qquad \zeta\sim\mathcal N(0,1),
\]
where $Z$ and $\zeta$ are independent. The representation follows from a
Brownian exponential functional and Dufresne's identity \cite{MY}; related
additive--multiplicative Cauchy statistics appear in \cite{DM}.

The limit law has Cauchy tails in every nontrivially forced direction. Thus a
global second-moment comparison is unavailable at stationarity. Instead, we
couple the nonlinear and linear equations with the same noises and compare them
up to the exit from a fixed $V$-ball. The stopped difference admits an $L^2$
estimate. The Lyapunov bound controls the exit probability, and the limiting
linear flow contracts fractional moments. Balancing these three estimates gives
a quantitative rate that is uniform over the whole invariant-measure family.

\subsection{Main results and organization}
We summarize the main conclusions. The emphasis is on estimates uniform over
all invariant probabilities, rather than on selecting a single stationary
measure: no uniqueness or ergodicity statement for the nonlinear dynamics is
needed for the small-forcing limit.

\begin{enumerate}[label=\textup{(\roman*)},leftmargin=*]
	
	\item \emph{Local and global well-posedness.}
	Under Assumption~\ref{ass:local-structure} and the regularization hypotheses
	(R1)--(R4), Theorem~\ref{thm:local-criterion} gives a unique maximal local
	solution with paths in $V$. Pathwise uniqueness and continuous dependence hold
	in $X$, and the high-order energy identity is recovered from the regularized
	equations. Under the radial growth condition above,
	Theorem~\ref{thm:global} gives global well-posedness and fractional Lyapunov
	bounds that are uniform in $\eta$. At $\eta=0$, the $q$th moment, where
	$q=2p$, decays exponentially, and $\delta_0$ is the unique invariant
	probability.
	
	\item \emph{Invariant measures and the stationary scale.}
	For every $\eta>0$, the Markov dynamics associated with \eqref{eq:original}
	admit invariant probabilities. Let $\mathcal I_\eta$ denote the full set of
	such probabilities, viewed as measures on $X$. Theorem~\ref{thm:invariant}
	proves the uniform estimate
	\begin{equation}\label{eq:intro-scale}
		\sup_{\mu\in\mathcal I_\eta}
		\int_V(1+\rho(u))(\eta^2+\|u\|_V^2)^p\,\mu(\dd u)
		\le C\eta^{2p}.
	\end{equation}
	Hence $\eta$ is the stationary scale for every invariant measure, not only for
	one measure selected by the Krylov--Bogoliubov procedure.
	Proposition~\ref{prop:stationary-set} shows that the invariant-measure sets are
	nonempty, convex, and compact. It also gives concentration at the origin as
	$\eta\downarrow0$.
	
	\item \emph{A universal rescaled Cauchy law.}
	Define
	\[
	D_\eta u=\frac{u}{\eta},
	\qquad
	\mathcal J_\eta=(D_\eta)_\#\mathcal I_\eta.
	\]
	If $G\ne0$, Theorem~\ref{thm:linear-limit} identifies $\nu_*$ as the unique
	invariant law of the limiting linear equation. Theorem~\ref{thm:scaling} then
	gives, for $0<\eta\le e^{-1}$,
	\begin{equation}\label{eq:intro-rate}
		\sup_{\nu\in\mathcal J_\eta}
		d_{\mathrm{BL},X}(\nu,\nu_*)
		\le K(1+\log\eta^{-1})\eta^\vartheta,
	\end{equation}
	where $d_{\mathrm{BL},X}$ is the bounded--Lipschitz metric and
	$\vartheta>0$ is explicit. The diameter of $\mathcal J_\eta$ therefore tends
	to zero, even if the invariant measure is not known to be unique at fixed
	$\eta$. For every $0<\alpha<q$, fractional $X$-moments and the Wasserstein
	distance associated with $\|x-y\|_X^\alpha$ also converge uniformly. The same
	coupling proves convergence of stationary path laws on every finite time
	interval.
	
	\item \emph{Applications and fixed-forcing tails.}
	Section~\ref{sec:models} verifies the abstract assumptions for active scalars,
	including two-dimensional Euler and SQG, two-dimensional ideal MHD,
	smoothing-inertia EPDiff equations, and the Camassa--Holm,
	Degasperis--Procesi, and Hunter--Saxton equations. In each case, the geometric
	formulation identifies the high-order cancellation and a conserved lower-order
	quadratic energy.
	
	When this lower-order energy is forced in at least one direction,
	Theorem~\ref{thm:low} yields, for every fixed $\eta>0$ and every
	$\mu\in\mathcal I_\eta$,
	\[
	\int_V \psi(\rho(u))^2\sqrt{1+E_0(u)}\,\mu(\dd u)=\infty.
	\]
	For square-root feedback in the models above, this implies
	$\int_V\|u\|_V^2\,\mu(\dd u)=\infty$. Thus the small-forcing limit is naturally
	formulated in fractional moments, even though the stationary high-order second
	moment diverges at each fixed nonzero forcing amplitude.
	
\end{enumerate}

From the viewpoint of stochastic dynamics, the results describe the
small-parameter behavior of the full family of stationary statistical states,
not merely the convergence of selected stationary solutions. In particular,
the limiting singleton provides a parameter-dependent selection principle at
the level of invariant probabilities, despite possible nonuniqueness of the
nonlinear stationary dynamics for every fixed positive forcing amplitude.

The central analytical contribution is a passage from global solvability
under radial stabilization to the stationary small-forcing regime for
inviscid derivative-losing dynamics with heavy-tailed stationary laws. In particular,
the uniform limit holds without assuming uniqueness, irreducibility, or
a strong Feller property for the nonlinear Markov dynamics at fixed
positive $\eta$. The proof combines
(i) a two-space regularization argument that recovers the high-order energy
identity, (ii) a truncated stationary Lyapunov argument that yields moment
bounds for every invariant probability without assuming its moments are finite
in advance, and (iii) a stopped common-noise comparison with the limiting linear stochastic equation.
The last step uses only fractional stationary moments, since a global
second-moment coupling is incompatible with the Cauchy tails of the limit.

This distinguishes the present results from the nonlinear-noise
non-explosion theory of Bagnara, Maurelli, and Xu \cite{BMX}, which addresses
global solvability but not the uniform small-forcing behavior of the full
invariant-measure set, and from the Cauchy stationary statistics for
additive--multiplicative stochastic heat equations studied by Dunlap and
Mukherjee \cite{DM}, which concern a dissipative parabolic model rather than
inviscid derivative-losing transport. The study of invariant probabilities in weak topologies has a long
history in stochastic evolution equations, including \cite{MS94,MS98}. For
a more recent Camassa--Holm setting involving continuity on regularity-bounded
sets, see \cite{APT}. Here we prove the specific bounded-$V$-ball,
$X$-topology continuity needed to pass through the transition kernels, and
combine it with a stationary estimate uniform in the forcing parameter.
The geometric structures of the applications supply the model-dependent
energy cancellations, whereas the limiting law depends only on
$a_0=\psi(0)$ and the additive covariance $Q_V=GG_V^*$.

	The paper is organized as follows. Section~\ref{sec:framework} develops the two-space local theory, recovers the
high-order energy identity, and proves global stabilization.
Section~\ref{sec:invariant} constructs invariant measures and establishes
stationary compactness at the forcing scale. Section~\ref{sec:scaling}
identifies the limiting linear law and proves convergence of invariant measures,
fractional moments, Wasserstein distances, and stationary paths.
Section~\ref{sec:models} verifies the hypotheses for the concrete equations.
Appendix~\ref{sec:low} proves the fixed-forcing moment obstruction, and
Appendix~\ref{sec:linear} treats a limiting equation in which a deterministic
linear mode remains after rescaling.

\medskip 
	
	\section{Well-posedness}\label{sec:framework}
	This section introduces the $V$-valued, $X$-stable solution concept and the
	analytic assumptions on the drift.  We construct maximal local solutions by
	regularization, recover the high-order energy identity, and then use the radial
	growth condition to prove global well-posedness.
	
	\subsection{Phase spaces and assumptions}
	We use two Hilbert spaces to separate spatial regularity from stability.  The
	space $V$ carries the high-order energy and the solution paths, while $X$ is used
	for difference estimates and weak convergence of probability measures.  Let
	$i:V\hookrightarrow X$ be a continuous, injective, compact embedding of real
	separable Hilbert spaces, and set
	\[
	K_R=\{u\in V:\|u\|_V\le R\},\qquad R>0.
	\]
	We identify $V$ with its image in $X$.  Each $K_R$ is compact and closed in $X$,
	and the extension of $\|\cdot\|_V$ by $+\infty$ outside $V$ is lower
	semicontinuous on $X$.  Moreover, $V=\bigcup_{m\ge1}K_m$ is Borel in $X$.
	The continuous injection $i$ between these Polish spaces is a Borel
	isomorphism onto its image, so the strong Borel sigma-field of $V$ agrees with
	the trace of $\mathcal B(X)$.
	
	We work on a complete filtered probability space
	$(\Omega,\mathcal F,(\mathcal F_t)_{t\ge0},\PP)$ satisfying the usual conditions.
	It carries an $(\mathcal F_t)$-Brownian motion $W$ and an independent cylindrical
	$(\mathcal F_t)$-Brownian motion $B$ on a separable Hilbert space $\mathfrak U$.
	Let $G\in\mathcal L_2(\mathfrak U,V)$ and define
	\[
	G_X=iG,\qquad
	Q_V=GG_V^*,\qquad Q_X=G_XG_X^*,\qquad
	q_V=\|G\|_{\mathcal L_2(\mathfrak U,V)}^2,
	\]
	where $G_V^*:V\to\mathfrak U$ and $G_X^*:X\to\mathfrak U$ denote the adjoints in
	the indicated Hilbert structures. We identify $X^*$ with $X$ through the
	Riesz map when applying $G_X^*$ to a
	continuous linear functional. Initial variables are $\mathcal F_0$-measurable
	and independent of future noise increments. Unless stated otherwise, constants
	may change from line to line and are independent of $\eta\in[0,1]$.
	
	\begin{definition}[$V$-valued, $X$-stable local solution]\label{def:vx-wp}
		For $u\in V$, a maximal local solution is a pair $(U^u,\tau^{*,u})$, where
		$\tau^{*,u}$ is a strictly positive stopping time and $U^u$ is progressively
		measurable with $U^u\in C([0,\tau^{*,u});V)$ almost surely, such that
		\begin{equation}\label{eq:vx-integral}
			U^u(t)=u-\int_0^tF(U^u(s))\dd s
			+\int_0^ta(U^u(s))U^u(s)\dd W_s
			+\eta\int_0^tG\dd B_s
		\end{equation}
		holds in $X$ for $t<\tau^{*,u}$, where $a(u)=\psi(\rho(u))$.  For $M<R$ define
		\[
		\tau_R^{u,v}=\inf\{t<\tau^{*,u}\wedge\tau^{*,v}:
		\|U^u(t)\|_V\vee\|U^v(t)\|_V\ge R\}.
		\]
		Local well-posedness in the $V$-valued, $X$-stable sense means pathwise
		uniqueness, Borel dependence of the stopped path laws on the initial state, the
		$V$-norm blow-up alternative, and
		\begin{equation}\label{eq:local-stability}
			\E\sup_{t\le T\wedge\tau_R^{u,v}}
			\|U^u(t)-U^v(t)\|_X^2
			\le C_{M,R,T}\|u-v\|_X^2,
			\qquad u,v\in K_M.
		\end{equation}
	\end{definition}
	
	\begin{assumption}[Local analytic structure]\label{ass:local-structure}
		The following conditions hold.
		\begin{enumerate}[label=\textup{(H\arabic*)},leftmargin=*]
			\item\label{hyp:rho}
			$\rho:V\to[0,\infty)$ is one-homogeneous and Lipschitz in the $X$ topology:
			\begin{equation}\label{eq:rho}
				\rho(\lambda u)=\lambda\rho(u)\quad(\lambda\ge0),\qquad
				|\rho(u)-\rho(v)|\le L_\rho\|u-v\|_X.
			\end{equation}
			In particular, $\rho(u)\le C_\rho\|u\|_V$.
			
			\item\label{hyp:drift}
			$F:V\to X$ is continuous and quadratic,
			\begin{equation}\label{eq:quadratic}
				F(\lambda u)=\lambda^2F(u),\qquad
				\|F(u)\|_X\le C_F\|u\|_V^2,
			\end{equation}
			and, for $u,v\in K_R$,
			\begin{equation}\label{eq:diff}
				\left|2\langle F(u)-F(v),u-v\rangle_X\right|
				\le C_R\|u-v\|_X^2.
			\end{equation}
			
			\item\label{hyp:energy-pairing}
			There is a dense smooth core $\mathscr D\subset V$ such that
			$F(\mathscr D)\subset V$ and
			\begin{equation}\label{eq:b-definition}
				b_V^0(u)=-2\langle F(u),u\rangle_V,
				\qquad u\in\mathscr D,
			\end{equation}
			extends to a continuous cubic functional $b_V:V\to\R$ satisfying
			\begin{equation}\label{eq:b}
				b_V(\lambda u)=\lambda^3b_V(u),\qquad
				|b_V(u)|\le C_V\rho(u)\|u\|_V^2.
			\end{equation}
			When $F(u)$ is only $X$-valued, \eqref{eq:b-definition} is understood through
			the regularized limit \eqref{eq:recover}. Condition \eqref{eq:b} is the
			model-specific top-order cancellation used in the high-order energy estimate.
			
			\item\label{hyp:psi-local}
			$\psi:[0,\infty)\to\R$ is locally Lipschitz, and
			$a(u)=\psi(\rho(u))$.
		\end{enumerate}
	\end{assumption}
	
	\begin{assumption}[Radial stabilization]\label{ass:stabilization}
		There exist $b_0,b_1>0$ such that
		\begin{equation}\label{eq:feedback}
			\psi(y)^2\ge b_0+b_1y,
			\qquad y\ge0.
		\end{equation}
		We take
		\begin{equation}\label{eq:a0}
			a_0=\psi(0)>0.
		\end{equation}
		The case $\psi(0)<0$ is reduced to this one by replacing $W$ with $-W$.
		Fix $p\in(0,1/2)$ such that
		\begin{equation}\label{eq:strength}
			(1-2p)b_1>C_V,
		\end{equation}
		and set
		\[
		\kappa=1-2p,
		\qquad \PhiV(u)=(1+\|u\|_V^2)^p,
		\qquad q=2p\in(0,1).
		\]
	\end{assumption}
	
	For $\eta>0$ define the dilation
	\begin{equation}\label{eq:dilation}
		D_\eta:X\to X,\qquad D_\eta u=\eta^{-1}u,
		\qquad (D_\eta)_\#\mu(A)=\mu(\eta A).
	\end{equation}
	The scaled variable $v_\eta=D_\eta U_\eta$ solves
	\begin{equation}\label{eq:scaled}
		\dd v_\eta+\eta F(v_\eta)\dd t
		=a_\eta(v_\eta)v_\eta\dd W_t+G\dd B_t,
		\qquad a_\eta(v)=\psi(\eta\rho(v)).
	\end{equation}
	At $\eta=0$ the limiting linear equation is
	\begin{equation}\label{eq:linear-limit}
		\dd V=a_0V\dd W_t+G\dd B_t.
	\end{equation}
	
	\subsection{Local well-posedness}\label{sec:local-criterion}
	We construct the local solution from regularized equations on the original
	stochastic basis.  The proof combines uniform high-order bounds, convergence in
	$X$, and recovery of the $V$-energy identity for the limit process.  Related
	regularization and lower-topology arguments appear in stochastic fluid theories
	\cite{GHV,CFH,ART,LC}; for Camassa--Holm equations, see also \cite{RT,MRT}.
	
	We specify the regularizers used in this section. Let
	$j\in C_c^\infty(\R^d;[0,1])$ be real and even, with $j=1$ near the origin,
	and define the componentwise Fourier multiplier
	\begin{equation}\label{eq:Jn-symbol}
		\widehat{J_nf}(k)=j(k/n)\widehat f(k),\qquad k\in\mathbb Z^d.
	\end{equation}
	We work with closed constraint subspaces preserved by these multipliers.
	The smooth core in condition~\ref{hyp:energy-pairing} is chosen to contain
	$J_nX$. In particular,
	\begin{align}\label{eq:Jn}
		J_n:X\longrightarrow\mathscr D\subset V,\qquad
		J_n^*=J_n\text{ in both }X\text{ and }V,\qquad
		[J_n,\Lambda^a]=0,\notag\\
		\|J_n\|_{\mathcal L(H^a)}\le1,\quad
		J_nf\longrightarrow f\text{ in }H^a,\quad
		\sup_n\|J_n\|_{\mathcal L(W^{1,\infty})}<\infty.
	\end{align}
	Here $\Lambda^a=(I-\Delta)^{a/2}$, and $J_n$ has finite-dimensional range.
	For all $b>a$,
	\begin{equation}\label{eq:mollifier}
		\|(I-J_n)f\|_{H^a}\le Cn^{-(b-a)}\|f\|_{H^b}.
	\end{equation}
	For fixed Lipschitz $v$ and $f\in L^2$, the periodic kernel estimate gives
	\begin{equation}\label{eq:mollifier-comm}
		\|[J_n,v]\partial_jf\|_{L^2}
		\le C\|v\|_{W^{1,\infty}}\|f\|_{L^2},\qquad
		[J_n,v]\partial_jf\longrightarrow0\text{ in }L^2.
	\end{equation}
	Indeed, with $K_n$ the convolution kernel,
	\[
	[J_n,v]\partial_jf(x)
	=\int_{\T^d}\partial_jK_n(x-y)(v(y)-v(x))f(y)\dd y
	-J_n((\partial_jv)f)(x).
	\]
	The uniform kernel bounds give the inequality. Convergence holds first for smooth $f$ and then for every $f\in L^2$ by
	density. Thus the convergence is strong for each fixed pair $v,f$.

	Let $V=H^s$ and $X=H^r$ on a torus, or finite products and Fourier-invariant
	closed linear constraint subspaces, with $r<s-1$.  Assume that the regularizers
	in \eqref{eq:Jn-symbol} satisfy
	\eqref{eq:Jn}--\eqref{eq:mollifier-comm}, preserve the phase space,
	and have the following properties:
	\begin{enumerate}[label=\textup{(R\arabic*)},leftmargin=*]
		\item The maps $F_n:V\to V$ are Lipschitz on every bounded $V$-ball, with
		constants allowed to depend on $n$, and
		\[
		|2\ip{F_n(v)}{v}_V|\le C\rho(v)\norm{v}_V^2
		\]
		uniformly in $n$.
		\item For some $\delta>0$ and every $M$,
		\[
		\sup_{v\in K_M}\norm{F_n(v)-F(v)}_X\le C_M n^{-\delta}.
		\]
		\item On each bounded $V$-ball, $F$ is continuous from the $X$ topology to
		$X$; local H\"older continuity is sufficient.
		\item For every $v\in V$,
		\begin{equation}\label{eq:recover}
			-2\ip{J_nF(v)}{J_nv}_V\longrightarrow b_V(v),\qquad
			|2\ip{J_nF(v)}{J_nv}_V|
			\le C\rho(v)\norm{v}_V^2.
		\end{equation}
	\end{enumerate}
	The pairings in (R4) are well defined because $J_n$ regularizes the
	$X$-valued drift.
	
	The first result establishes local well-posedness in the strong–weak topology framework.
	
	\begin{theorem}[Local well-posedness]\label{thm:local-criterion}
		Under Assumption~\ref{ass:local-structure} and (R1)--(R4),
		\eqref{eq:original} is locally well posed in the $V$-valued, $X$-stable sense of
		Definition~\ref{def:vx-wp}.  The estimate \eqref{eq:local-stability} holds, and
		the localized path laws form Borel kernels on $V$. 
		If
		\[
		\sigma_R^u=\inf\{t<\tau^{*,u}:\rho(U^u(t))\ge R\},
		\]
		then $\sigma_R^u\uparrow\tau^{*,u}$ almost surely, and
		\begin{equation}\label{eq:rho-continuation}
			\{\tau^{*,u}<\infty\}
			\subset
			\left\{\limsup_{t\uparrow\tau^{*,u}}\rho(U^u(t))=\infty\right\}
			\subset
			\left\{\limsup_{t\uparrow\tau^{*,u}}\norm{U^u(t)}_V=\infty\right\}.
		\end{equation}
		The constants in the fixed-ball estimates can be chosen uniformly for the
		scaled coefficients $\eta F$ and $v\mapsto\psi(\eta\rho(v))v$, with
		$0\le\eta\le1$.
	\end{theorem}
	\begin{proof}
		\emph{Truncated Hilbert-space equations.}
		Let $\chi_R$ be a smooth cutoff equal to one on $[0,R]$ and zero on $[2R,\infty)$. Set
		\begin{align*}
			F_{n,R}(v)&=\chi_R(\rho(v))^2F_n(v),\\
			H_R(v)(z,e)&=\chi_R(\rho(v))\,[z\psi(\rho(v))v+\eta Ge],
			\qquad(z,e)\in\R\oplus\mathfrak U.
		\end{align*}
		For fixed $n$ and $R$, both maps are locally Lipschitz in $V$. Standard local
		Hilbert-space SDE theory \cite{DPZ} gives a solution of
		\begin{equation}\label{eq:reg-local}
			\dd u_n+F_{n,R}(u_n)\dd t=H_R(u_n)\dd(W,B),\qquad u_n(0)=u_0\in V.
		\end{equation}
		The same cutoff is applied to the additive term, so all coefficients of the
		approximating equation satisfy uniform energy bounds on the active region. The
		first regularization condition and local boundedness of $\psi$ give
		\[
		|2\ip{F_{n,R}(v)}{v}_V|\le C_R\norm{v}_V^2,\qquad
		\norm{H_R(v)}_{\HS(\R\oplus\mathfrak U,V)}^2
		\le C_R(1+\norm{v}_V^2).
		\]
		Apply It\^o's formula to $\|u_n\|_V^2$ and $\|u_n\|_V^4$, stopped at the
		norm-localization time. The martingale maximal inequality and Gronwall's lemma
		then give
		\begin{equation}\label{eq:local-moments}
			\sup_n\E\sup_{t\le T}\norm{u_n(t)}_V^{2j}
			\le C_{R,T,j}(1+\norm{u_0}_V^{2j}),\qquad j=1,2.
		\end{equation}
		These bounds preclude explosion of each regularized cutoff equation.
		
		\emph{Cauchy convergence in the lower space.}
		Stop $u_n,u_m$ at
		\[
		\tau_M^{n,m}=\inf\{t\ge0:\|u_n(t)\|_V\vee\|u_m(t)\|_V\ge M\}.
		\]
		The drift can be written as
		$\chi_R(\rho(v))^2[F(v)+(F_n(v)-F(v))]$.
		Set $w_{n,m}=u_n-u_m$.  Using \eqref{eq:diff}, the Lipschitz bound for the cutoff, and (R2), the localized $X$-energy satisfies
		\begin{align}
			\dd\norm{w_{n,m}}_X^2
			&\le C_{R,M}\norm{w_{n,m}}_X^2\dd t
			+C_M(n^{-2\delta}+m^{-2\delta})\dd t+\dd M_t,
			\label{eq:local-difference-ito}\\
			\dd[M]_t
			&\le 4\norm{w_{n,m}}_X^2
			\norm{H_R(u_n)-H_R(u_m)}_{\HS(\R\oplus\mathfrak U,X)}^2\dd t
			\le C_{R,M}\norm{w_{n,m}}_X^4\dd t.
			\label{eq:local-difference-bracket}
		\end{align}
		Here Young's inequality controls the approximation remainders, while the local
		Lipschitz bound for $H_R$ gives the bracket estimate.  With
		$\tau=\tau_M^{n,m}$, the Burkholder--Davis--Gundy (BDG) inequality yields
		\begin{align*}
			\E\sup_{t\le T}|M_{t\wedge\tau}|
			&\le C\E\left[
			\sup_{t\le T\wedge\tau}\|w_{n,m}(t)\|_X^2
			\int_0^{T\wedge\tau}\|w_{n,m}(s)\|_X^2\dd s
			\right]^{1/2}\\
			&\le \frac12\E\sup_{t\le T\wedge\tau}\|w_{n,m}(t)\|_X^2
			+C_{R,M}\int_0^T
			\E\sup_{r\le s\wedge\tau}\|w_{n,m}(r)\|_X^2\dd s.
		\end{align*}
		Consequently, Gronwall's lemma gives
		\begin{equation}\label{eq:local-cauchy}
			\E\sup_{t\le T\wedge\tau_M^{n,m}}
			\norm{u_n(t)-u_m(t)}_X^2
			\le C_{R,M,T}(n^{-2\delta}+m^{-2\delta}).
		\end{equation}
		For every $\varepsilon>0$, \eqref{eq:local-moments} and
		\eqref{eq:local-cauchy} give
		\begin{align*}
			\PP\!\left(\sup_{t\le T}\|u_n(t)-u_m(t)\|_X>\varepsilon\right)
			&\le \PP(\tau_M^{n,m}\le T)
			+\varepsilon^{-2}\E\sup_{t\le T\wedge\tau_M^{n,m}}
			\|u_n(t)-u_m(t)\|_X^2\\
			&\le \frac{C_{R,T}(1+\|u_0\|_V^2)}{M^2}
			+\frac{C_{R,M,T}}{\varepsilon^2}
			(n^{-2\delta}+m^{-2\delta}).
		\end{align*}
		Letting first $n,m\to\infty$ and then $M\to\infty$ shows that $u_n$ is
		Cauchy in probability in $C([0,T];X)$ on the original stochastic basis.
		
		\emph{Identification and top-order continuity.}
		Choose a subsequence that converges almost surely in $C([0,T];X)$. To recover
		the top-order bound, extend $\norm{\cdot}_V$ by $+\infty$ outside $V$ and use
		the path functional
		\[
		\mathcal M_T(w)=\sup_{t\in([0,T]\cap\mathbb Q)\cup\{T\}}
		\norm{w(t)}_V.
		\]
		The functional $\mathcal M_T$ is lower semicontinuous on $C([0,T];X)$.
		Fatou's lemma and \eqref{eq:local-moments} show that $\mathcal M_T(u)$ is
		almost surely finite and satisfies the same second- and fourth-moment bounds.
		Since each $K_M$ is closed in $X$, continuity in $X$ extends the bound from
		rational times to all times. On a bounded $V$-ball, convergence in $X$ implies
		weak convergence in $V$ by compactness and uniqueness of the embedded limit.
		Thus $u$ is weakly continuous in $V$; the Borel identification above gives
		strong measurability.
		
		Let $F_R(v)=\chi_R(\rho(v))^2F(v)$. By (R3), $F_R|_{K_M}$ is uniformly
		continuous in $X$, since $K_M$ is compact. Let $\omega_{R,M}$ be its modulus
		of continuity, with $\omega_{R,M}(0+)=0$, and set
		\[
		A_{n,M}=\left\{\sup_{t\le T}
		(\|u_n(t)\|_V\vee\|u(t)\|_V)\le M\right\}.
		\]
		The moment bounds and Fatou's lemma imply
		$\sup_n\PP(A_{n,M}^c)\le C_{R,T}(1+\|u_0\|_V^2)M^{-2}$. On $A_{n,M}$,
		\begin{equation}\label{eq:localized-drift-identification}
			\sup_{t\le T}\left\|\int_0^t
			[F_{n,R}(u_n)-F_R(u)]\dd s\right\|_X
			\le C_MTn^{-\delta}
			+T\omega_{R,M}\!\left(\sup_{t\le T}\|u_n-u\|_X\right).
		\end{equation}
		Sending $n\to\infty$ and then $M\to\infty$ proves convergence of the drift
		integrals in probability.
		
		On $K_M$, $H_R$ is Lipschitz into $\HS(\R\oplus\mathfrak U,X)$.
		Consequently $H_R(u_n)\to H_R(u)$ in probability in
		$L^2(0,T;\HS(\R\oplus\mathfrak U,X))$. The estimate
		\[
		\int_0^T\|H_R(u_n)-H_R(u)\|_{\HS}^2\dd s
		\le C_{R,T}\left(1+\sup_{t\le T}\|u_n(t)\|_V^2
		+\sup_{t\le T}\|u(t)\|_V^2\right)
		\]
		and the fourth-moment bound give uniform integrability. The convergence
		therefore holds in expectation, and the maximal inequality gives convergence of
		the stochastic integrals in $L^2(\Omega;C([0,T];X))$. Passing to the limit
		yields the truncated equation on the original stochastic basis.
		
		Put $Y_k=\|J_ku\|_V^2$, $\gamma_R=\chi_R(\rho(u))$, and
		$q_{V,k}=\|J_kG\|_{\HS(\mathfrak U,V)}^2$.  Applying It\^o's formula to
		$J_ku$ gives
		\begin{align}
			\dd Y_k={}&\gamma_R^2\!
			\left[-2\ip{J_kF(u)}{J_ku}_V+a(u)^2Y_k+\eta^2q_{V,k}\right]\dd t\notag\\
			&+2\gamma_Ra(u)Y_k\dd W_t
			+2\eta\gamma_R\ip{(J_kG)_V^*J_ku}{\dd B_t}_{\mathfrak U}.
			\label{eq:local-recovered-ito}
		\end{align}
		Write $Y=\|u\|_V^2$ and
		\[
		B_k(u)=-2\ip{J_kF(u)}{J_ku}_V,\qquad
		h_k=(J_kG)_V^*J_ku=G_V^*J_k^2u.
		\]
		On the active region of the cutoff, (R4) and local boundedness of $\psi$ give
		\[
		\gamma_R^2(|B_k(u)|+|b_V(u)|)\le C_RY,\qquad
		|\gamma_R a(u)|\le C_R,\qquad |Y_k-Y|\le2Y.
		\]
		Since $\E\sup_{t\le T}Y(t)^2<\infty$, dominated convergence yields
		\begin{align}
			&\E\int_0^T\gamma_R^2|B_k(u)-b_V(u)|\dd s\longrightarrow0,
			\label{eq:recovery-drift-L1}\\
			&\E\int_0^T\gamma_R^2a(u)^2|Y_k-Y|^2\dd s\longrightarrow0,
			\qquad q_{V,k}\longrightarrow q_V.
			\label{eq:recovery-radial-L2}
		\end{align}
		Moreover,
		\[
		\|h_k-G_V^*u\|_{\mathfrak U}
		\le \|G\|_{\HS(\mathfrak U,V)}\|(J_k^2-I)u\|_V
		\longrightarrow0,
		\]
		and the square of the left-hand side is bounded by $4q_VY$. If $M_k$ and $M$
		denote the sums of the two martingales in the regularized and limiting energy
		formulas, Doob's inequality
		and independence of $W,B$ imply
		\begin{align}\label{eq:recovery-martingale-L2}
			\E\sup_{t\le T}|M_k(t)-M(t)|^2
			\le{}&16\E\int_0^T\gamma_R^2a(u)^2|Y_k-Y|^2\dd s\notag\\
			&+16\eta^2\E\int_0^T\gamma_R^2
			\|h_k-G_V^*u\|_{\mathfrak U}^2\dd s
			\longrightarrow0.
		\end{align}
		The remaining finite-variation terms converge in $L^1(\Omega;C([0,T]))$ by
		the same domination. Thus the integrated right-hand sides of
		\eqref{eq:local-recovered-ito} converge uniformly in time in probability to a
		continuous process $\mathcal Y$. Choose a subsequence for which this convergence
		is almost surely uniform. On the event where the path is bounded in $V$,
		$J_ku(t)\to u(t)$ in $V$ at every $t$, so
		\[
		\mathcal Y(t)=\|u(t)\|_V^2\quad(0\le t\le T).
		\]
		Weak $V$ continuity and norm continuity imply $u\in C([0,T];V)$. The limiting
		energy martingale has integrands $2\gamma_Ra(u)Y$ and
		$2\eta\gamma_RG_V^*u$. Its quadratic variation follows from their squares and
		the zero cross-variation of $W$ and $B$. Before the feedback exit,
		$\gamma_R=1$, and the limit satisfies
		\eqref{eq:energy}--\eqref{eq:bracket}.
		
		\emph{Stability with respect to initial states.}
		Let $u^x,u^y$ be two limiting solutions driven by the same noises and
		started from $x,y\in K_M$.  Stop them at $\tau_R^{x,y}$ with $R>M$.  On this
		interval the additive noises cancel, while the drift difference obeys
		\eqref{eq:diff} and the map $u\mapsto a(u)u$ is Lipschitz in $X$ on $K_R$.
		It\^o's formula gives, after localization,
		\[
		\dd\norm{u^x-u^y}_X^2
		\le C_R\norm{u^x-u^y}_X^2\dd t
		+2\ip{u^x-u^y}{a(u^x)u^x-a(u^y)u^y}_X\dd W_t.
		\]
		The BDG inequality, Young's inequality, and Gronwall's
		lemma yield
		\begin{equation}\label{eq:criterion-stability}
			\E\sup_{t\le T\wedge\tau_R^{x,y}}
			\norm{u^x(t)-u^y(t)}_X^2
			\le C_{M,R,T}\norm{x-y}_X^2,
		\end{equation}
		which is \eqref{eq:local-stability}. The same estimate gives pathwise
		uniqueness and continuous dependence in the lower topology.
		
		\emph{Uniqueness, maximality, continuation, and measurability.}
		The stopped $X$-difference estimate proves pathwise uniqueness. For $R>\rho(u_0)$, the truncated solution solves the original equation until
		$\rho(u)$ reaches $R$.  The cutoff equation is global: on its active region
		the high-order energy has at most linear growth, and outside that region
		all truncated stochastic coefficients vanish.  Patching these global
		cutoff solutions by uniqueness gives a maximal solution and shows that its
		lifetime is the increasing limit of the feedback exit times $\sigma_R^u$.
		If that lifetime were finite while $\rho(u)$ remained bounded, a larger
		cutoff solution would extend it.  This proves the first inclusion in
		\eqref{eq:rho-continuation}; the second follows from
		$\rho(u)\le C_\rho\norm{u}_V$.
		
		\emph{Borel path-law kernels and random initial data.}
		Fix $R,T$, and denote by $u_R^x$ the global cutoff solution. Repeating the
		localized difference estimate for this equation gives, for $x,y\in K_L$ and
		$M>L$,
		\begin{equation}\label{eq:cutoff-initial-kernel-bound}
			\E\left[1\wedge\|u_R^x-u_R^y\|_{C([0,T];X)}^2\right]
			\le C_{R,M,T}\|x-y\|_X^2
			+\frac{C_{R,T}(1+L^2)}{M^2}.
		\end{equation}
		The cutoff in the additive term contributes an additional difference, which is
		controlled by the Lipschitz property of $\rho$; the exit bound follows from
		\eqref{eq:local-moments}.
		If $x_j\to x$ in $V$, the initial norms are bounded. Letting $j\to\infty$
		and then $M\to\infty$ shows that
		\begin{equation}\label{eq:path-law-kernel}
			x\longmapsto\mathcal K_{R,T}(x,\cdot)
			:=\Law(u_R^x|_{[0,T]})
			\quad\text{is continuous from }V\text{ to }
			\mathcal P(C([0,T];X)).
		\end{equation}
		It is therefore a Borel probability kernel. The injection
		$C([0,T];V)\to C([0,T];X)$ is continuous and injective between Polish spaces,
		so its image and inverse are Borel. Consequently, stopping the cutoff paths at a feedback exit also gives Borel
		localized path-law kernels.
		
		For a bounded $V$-valued $\mathcal F_0$-measurable variable $\xi$, choose
		countably valued $\xi_j$ with $\|\xi_j-\xi\|_V\le j^{-1}$ and bounded $V$
		norms. Construct $u_R^{\xi_j}$ by patching the deterministic-initial solutions
		on the $\mathcal F_0$ events defining $\xi_j$. The estimates above, the
		identification argument, and \eqref{eq:recovery-drift-L1}--\eqref{eq:recovery-martingale-L2}
		construct $u_R^\xi$ on the given stochastic basis. For $F\in C_b(C([0,T];X))$,
		independence of future noise increments first gives, for countably valued data,
		\[
		\E[F(u_R^{\xi_j})\mid\mathcal F_0]
		=\mathcal K_{R,T}F(\xi_j).
		\]
		Convergence in probability, boundedness, and \eqref{eq:path-law-kernel} pass this
		identity to $\xi$. A monotone-class argument extends it to bounded Borel $F$.
		Patching over the disjoint events $\{j-1\le\|\xi\|_V<j\}$ extends the
		construction to every almost surely $V$-valued initial variable $\xi$.
		
		Passing from the cutoff equations to the original equation gives the same
		properties for the maximal local solution. Once global existence is established,
		the cutoff path laws converge, for each initial point, to the global path law;
		the pointwise weak limit is again a Borel kernel. Restarting at
		a deterministic time $s$ with initial variable $U(s)$ and the independent shifted
		noises, pathwise uniqueness and the preceding conditional identity give
		\begin{equation}\label{eq:markov-restart}
			\E[f(U(s+t))\mid\mathcal F_s]=\widehat P_t^\eta f(U(s)),
			\qquad f\in B_b(X),\quad s,t\ge0.
		\end{equation}
		The global solution therefore has time-homogeneous Borel transition kernels on
		$V$, viewed as a Borel subset of $X$. For the scaled coefficients $\eta F$ and
		$\psi(\eta\rho)$, the fixed-ball bounds are uniform in
		$\eta\in[0,1]$.
	\end{proof}
	
	For $Y=\|U_\eta\|_V^2$, the solution constructed above satisfies, up to its
	lifetime,
	\begin{align}
		\dd Y={}&[b_V(U_\eta)+a(U_\eta)^2Y+\eta^2q_V]\dd t
		+2a(U_\eta)Y\dd W_t
		+2\eta\langle G_V^*U_\eta,\dd B_t\rangle_{\mathfrak U},\label{eq:energy}\\
		\dd[Y]_t={}&4\big[a(U_\eta)^2Y^2+
		\eta^2\|G_V^*U_\eta\|_{\mathfrak U}^2\big]\dd t.\label{eq:bracket}
	\end{align}
	
	\subsection{Global well-posedness}\label{sec:globalization}
	The local construction above identifies the only high-order term that enters the
	Lyapunov calculation, namely $b_V$.  The following theorem shows that the concavity
	of $(1+\|u\|_V^2)^p$ converts the radial quadratic variation into a coercive
	drift. Let $\mathscr L_\eta$ denote the generator of the scaled equation
	\eqref{eq:scaled}.
	
	\begin{theorem}[Global well-posedness and Lyapunov bounds]\label{thm:global}
		Let $(U_\eta,\tau_\eta^*)$ be the maximal local solution supplied by
		Theorem~\ref{thm:local-criterion}. Under
		Assumptions~\ref{ass:local-structure} and~\ref{ass:stabilization}, this
		solution is global.
		There are constants $c,C>0$, independent of $\eta\in[0,1]$, such that
		\begin{equation}\label{eq:lyap}
			\mathscr L_\eta\PhiV(v)
			\le C-c(1+\eta\rho(v))\PhiV(v)
			\le C-c\PhiV(v).
		\end{equation}
		For \eqref{eq:original}, the same constants give the corresponding bound with
		$1+\rho(v)$ in place of $1+\eta\rho(v)$, uniformly for
		$0\le\eta\le1$.
		All equations \eqref{eq:scaled}, including $\eta=0$, satisfy
		\begin{equation}\label{eq:moment}
			\E\PhiV(v_\eta(t))
			\le e^{-ct}\E\PhiV(v_\eta(0))+\frac Cc(1-e^{-ct})
		\end{equation}
		whenever the initial expectation is finite. If
		$\tau_R=\inf\{t\ge0:\norm{v_\eta(t)}_V\ge R\}$, with $\tau_R=0$ for an initial value outside $K_R$, then
		\begin{equation}\label{eq:exit}
			\PP(\tau_R\le T)
			\le\frac{\E\PhiV(v_\eta(0))+CT}{(1+R^2)^p}.
		\end{equation}
		Both estimates are uniform in $\eta$.
		
		If the additive amplitude in the original equation is zero, then, with $q=2p$,
		\begin{equation}\label{eq:zero-decay}
			\E\|U_0(t)\|_V^q
			\le e^{-p(1-2p)b_0t}\E\|U_0(0)\|_V^q,
		\end{equation}
		and $\delta_0$ is its unique invariant probability.
	\end{theorem}
	\begin{proof}
		\emph{Step 1: the scale-uniform Lyapunov drift.}
		Put $x=\norm{v}_V^2$. In the scaled equation the high-order drift is $\eta b_V(v)$ and the additive trace is $q_V$. It\^o's formula gives
		\begin{align}\label{eq:ito-phi}
			\mathscr L_\eta\PhiV(v)
			={}&p(1+x)^{p-1}[\eta b_V(v)+a_\eta(v)^2x+q_V]\notag\\
			&+2p(p-1)(1+x)^{p-2}
			[a_\eta(v)^2x^2+\norm{G_V^*v}_{\mathfrak U}^2].
		\end{align}
		The additive second-derivative term is nonpositive.  The two radial terms
		combine as
		\begin{align*}
			&p(1+x)^{p-1}a_\eta(v)^2x
			+2p(p-1)(1+x)^{p-2}a_\eta(v)^2x^2\\
			&\qquad=px(1+x)^{p-1}a_\eta(v)^2
			\frac{1-(1-2p)x}{1+x}.
		\end{align*}
		Using \eqref{eq:b} and $\kappa=1-2p$ therefore gives
		\begin{align}\label{eq:central}
			\mathscr L_\eta\PhiV(v)
			\le{}&px(1+x)^{p-1}
			\left[\eta C_V\rho(v)
			+\psi(\eta\rho(v))^2\frac{1-\kappa x}{1+x}\right]
			+p q_V(1+x)^{p-1}.
		\end{align}
		Choose $\delta$ with $C_V/b_1<\delta<\kappa$, and
		\[
		x_*\ge\max\left\{1,\frac{1+\delta}{\kappa-\delta}\right\}.
		\]
		For $x\ge x_*$, we have $(1-\kappa x)/(1+x)\le-\delta$, and hence
		\begin{align}\label{eq:negative-bracket}
			\eta C_V\rho(v)
			+\psi(\eta\rho(v))^2\frac{1-\kappa x}{1+x}
			&\le\eta C_V\rho(v)-\delta\psi(\eta\rho(v))^2\notag\\
			&\le-\delta b_0-(\delta b_1-C_V)\eta\rho(v).
		\end{align}
		Since $x/(1+x)\ge1/2$, this is a negative multiple of
		$(1+\eta\rho(v))\PhiV(v)$ after multiplication by the prefactor.
		For $x\le x_*$, the quantity $\eta\rho(v)$ is uniformly bounded by
		$C_\rho\sqrt{x_*}$. Local boundedness of $\psi$ therefore controls all
		positive terms. After increasing $C$, we obtain \eqref{eq:lyap}. Repeating the
		calculation for \eqref{eq:original} removes the factors $\eta$ from the drift
		and feedback arguments, while the additive trace satisfies
		$\eta^2q_V\le q_V$.
		
		\emph{Step 2: exits, non-explosion, and moments.}
		For a deterministic initial value, stop the energy formula at $\tau_R$ and
		localize the martingale. On the stopped $V$-ball, the drift and bracket rates are
		bounded, so the martingale has mean zero. The inequality $\mathscr L_\eta\PhiV\le C$ gives
		\[
		\E\PhiV(v_\eta(T\wedge\tau_R))\le\PhiV(v_\eta(0))+CT.
		\]
		Continuity yields \eqref{eq:exit}. The local blow-up alternative and $R\to\infty$ imply non-explosion. For the time-uniform bound, use the integrating factor \emph{before} removing the stopping:
		\[
		\E\left[e^{c(t\wedge\tau_R)}\PhiV(v_\eta(t\wedge\tau_R))\right]
		\le\PhiV(v_\eta(0))+\frac Cc(e^{ct}-1).
		\]
		Fatou's lemma gives \eqref{eq:moment}. Conditioning on the initial sigma-field
		extends the estimate to integrable random initial data; patching on bounded-norm
		events covers every almost surely $V$-valued initial variable.

		\emph{Step 3: the unforced endpoint.}
		For $Y=\|U_0\|_V^2$ and $q=2p$,
		\[
		\mathscr L_0^{\mathrm{orig}}Y^p
		=pY^{p-1}b_V(U_0)-p(1-2p)a(U_0)^2Y^p
		\le-p(1-2p)b_0Y^p.
		\]
		Stopping between the levels $\varepsilon$ and $R$, applying an integrating factor,
		and then sending $R\to\infty$ and $\varepsilon\downarrow0$ gives
		\eqref{eq:zero-decay}.  Hence $U_0^u(t)\to0$ in probability for every $u\in V$.
		If $\mu$ is invariant, then for $f\in C_b(V)$,
		\[
		\int f\,\dd\mu=\int P_t^{\mathrm{orig},0}f\,\dd\mu\longrightarrow f(0),
		\]
		so $\mu=\delta_0$.
	\end{proof}
	
\medskip 
	
	\section{Invariant measures and the stationary scale}\label{sec:invariant}
	We construct invariant measures from the global Lyapunov estimate and derive a
	small-forcing bound that is uniform over the entire invariant-measure family. The
	Markov dynamics start in $V$, while compactness and weak convergence are measured
	in $X$. Instead of assuming a global Feller property in $X$, we establish
	continuity of transition probabilities on bounded $V$-balls with the
	$X$-topology. This is related to the regularity-bounded continuity methods
	used for stochastic Camassa--Holm equations \cite{APT}; the precise
	continuity statement required here is proved below.
	
	\subsection{Invariant measures}
	For $f\in\mathcal B_b(X)$ and $u,v\in V$, let
	\begin{equation}\label{eq:explicit-transition-semigroups}
	\widehat P_t^\eta f(u)=\E[f(U_\eta^u(t))],\qquad
	P_t^\eta f(v)=\E[f(v_\eta^v(t))],
	\end{equation}
	where $U_\eta^u$ solves \eqref{eq:original} with initial state $u$,
	and $v_\eta^v$ solves \eqref{eq:scaled} with initial state $v$.
	The associated Borel transition kernels on $X$ are
	\begin{align}\label{eq:explicit-transition-kernels}
	\widehat P_t^\eta(u,A)&=\PP(U_\eta^u(t)\in A),\notag\\
	P_t^\eta(v,A)&=\PP(v_\eta^v(t)\in A),
	\qquad A\in\mathcal B(X).
	\end{align}
	The kernels are supported in $V$ when the initial condition belongs to $V$.
	We regard $V$ as a Borel subset of $X$ and consider only probabilities
	supported in $V$. For $\eta\in(0,1]$, define
	\begin{align}\label{eq:invariant-families}
		\mathcal I_\eta
		&=\{\mu\in\mathcal P(X):\mu(V)=1,
		\ \mu\widehat P_t^\eta=\mu\text{ for all }t\ge0\},\\
		\mathcal J_\eta
		&=\{\nu\in\mathcal P(X):\nu(V)=1,
		\ \nu P_t^\eta=\nu\text{ for all }t\ge0\}.
	\end{align}
	For a transition kernel $P_t$ and a probability $\mu$ supported in $V$, we use
	$\mu P_t(A)=\int_VP_t(u,A)\,\mu(\dd u)$.  For $\eta>0$,
	\begin{equation}\label{eq:invariant-dilation}
		\mathcal J_\eta=(D_\eta)_\#\mathcal I_\eta.
	\end{equation}
	
	The following theorem is the main stationary estimate: it establishes existence and a bound
	uniform over every invariant probability.
	
	\begin{theorem}[Invariant measures and stationary bounds]\label{thm:invariant}
	Assume the hypotheses of Theorem~\ref{thm:global}.
	Let $v_\eta^u$ denote the global solution of the scaled equation
	\eqref{eq:scaled} with $v_\eta^u(0)=u$, driven by the same noises $W$ and $B$.
	For fixed $T\ge0$, if $u_n,u\in K_M$ and $u_n\to u$ in $X$, then
	\begin{equation}\label{eq:global-restricted-path}
		\sup_{t\le T}\norm{v_\eta^{u_n}(t)-v_\eta^u(t)}_X
		\longrightarrow0\quad\text{in probability}.
	\end{equation}
	Consequently, $P_t^\eta f$ restricted to each $K_M$ is $X$-continuous for $f\in C_b(X)$.
	For every $\eta\in(0,1]$, the set $\cI_\eta$ is nonempty, and
	\begin{equation}\label{eq:scaled-stationary}
		\sup_{\mu\in\cI_\eta}
		\int_V(\eta^2+\norm{u}_V^2)^p\,\mu(\dd u)
		\le C_{\mathrm{st}}\eta^{2p},\qquad C_{\mathrm{st}}=C/c.
	\end{equation}
	More precisely, every $\mu\in\cI_\eta$ satisfies
	\[
	\int_V(1+\rho(u))(\eta^2+\norm{u}_V^2)^p\,\mu(\dd u)
	\le C_{\mathrm{st}}\eta^{2p}.
	\]
	Both estimates are uniform over $\mu\in\cI_\eta$.
\end{theorem}
\begin{proof}
	\emph{Step 1: restricted path continuity.}
	Write $H_\eta(v)=a_\eta(v)v$.  On $K_R$, uniformly in $\eta\in[0,1]$,
	\begin{align*}
		\|H_\eta(v)-H_\eta(w)\|_X
		&\le |a_\eta(v)|\|v-w\|_X
		+|a_\eta(v)-a_\eta(w)|\|w\|_X\\
		&\le L_R\|v-w\|_X.
	\end{align*}
	For $u_n,u\in K_M$ and $R>M$, put
	$z_n=v_\eta^{u_n}-v_\eta^u$ and let $\tau_R$ be the first time either process
	leaves $K_R$.  It\^o's formula and
	\eqref{eq:diff} imply, before $\tau_R$,
	\begin{align*}
		\dd\|z_n\|_X^2
		&=\left[-2\eta\ip{F(v_\eta^{u_n})-F(v_\eta^u)}{z_n}_X
		+\|H_\eta(v_\eta^{u_n})-H_\eta(v_\eta^u)\|_X^2\right]\dd t
		+\dd M_n\\
		&\le C_R\|z_n\|_X^2\dd t+\dd M_n,
		\qquad \dd[M_n]_t\le4L_R^2\|z_n\|_X^4\dd t.
	\end{align*}
	The BDG inequality applied to the stopped martingale, followed by
	Gronwall's inequality as in \eqref{eq:criterion-stability}, gives
	\[
	\E\sup_{t\le T\wedge\tau_R}\|z_n(t)\|_X^2
	\le C_{M,R,T}\|u_n-u\|_X^2.
	\]
	For any $\varepsilon>0$, the stopped estimate and \eqref{eq:exit} imply
	\begin{align}\label{eq:restricted-continuity-bound}
		\PP\!\left(\sup_{t\le T}\|z_n(t)\|_X>\varepsilon\right)
		&\le\varepsilon^{-2}\E\sup_{t\le T\wedge\tau_R}\|z_n(t)\|_X^2
		+\PP(\tau_R\le T)\notag\\
		&\le\frac{C_{M,R,T}}{\varepsilon^2}\|u_n-u\|_X^2
		+\frac{2[(1+M^2)^p+CT]}{(1+R^2)^p}.
	\end{align}
	Taking $n\to\infty$ and then $R\to\infty$ proves the
	convergence in probability in \eqref{eq:global-restricted-path}.
	For $f\in C_b(X)$, convergence in
	probability at time $t$ and boundedness give
	$P_t^\eta f(u_n)\to P_t^\eta f(u)$.
	
	\emph{Step 2: compactness and passage through the transition kernel.}
	Let $\lambda_n\in\mathcal P(X)$ be supported in $V$ and satisfy
	$\sup_n\int_V\PhiV\dd\lambda_n\le B$.  Then
	\begin{equation}\label{eq:compact-tail-bound}
		\lambda_n(K_R^c)\le\frac{B}{(1+R^2)^p}.
	\end{equation}
	Since $K_R$ is compact in $X$, this proves tightness.  If
	$\lambda_n\Rightarrow\lambda$ in $X$, extend $\PhiV$ by $+\infty$ outside
	$V$.  Lower semicontinuity and the Portmanteau theorem give
	\[
	\int_X\PhiV\dd\lambda
	\le\liminf_n\int_X\PhiV\dd\lambda_n\le B,
	\qquad \lambda(V)=1.
	\]
	Suppose $g:V\to\R$ is bounded and $g|_{K_R}$ is $X$-continuous.
	Tietze's extension theorem gives $g_R\in C_b(X)$ such that
	$g_R=g$ on $K_R$ and $\|g_R\|_\infty\le\|g\|_\infty$.  Therefore
	\begin{align}\label{eq:restricted-kernel-passage}
		\left|\int_Vg\dd\lambda_n-\int_Vg\dd\lambda\right|
		&\le\left|\int_Xg_R\dd\lambda_n-\int_Xg_R\dd\lambda\right|
		+2\|g\|_\infty[\lambda_n(K_R^c)+\lambda(K_R^c)]\notag\\
		&\le\left|\int_Xg_R\dd\lambda_n-\int_Xg_R\dd\lambda\right|
		+\frac{4B\|g\|_\infty}{(1+R^2)^p}.
	\end{align}
	Again let $n\to\infty$ and then $R\to\infty$. By Step~1, the same passage
	applies to $g=P_h^\eta f$.
	
	\emph{Step 3: construction of an invariant law.}
	Fix $v\in V$ and define the occupation measures
	\[
	\nu_T=\frac1T\int_0^T(P_t^\eta)^*\delta_v\dd t.
	\]
	The time-uniform estimate \eqref{eq:moment} gives
	\begin{equation}\label{eq:occupation-moment}
		\int_V\PhiV\dd\nu_T
		=\frac1T\int_0^T\E\PhiV(v_\eta^v(t))\dd t
		\le\frac{1-e^{-cT}}{cT}\PhiV(v)+\frac Cc.
	\end{equation}
	Step~2 yields $T_j\to\infty$ and $\nu_{T_j}\Rightarrow\nu$, with
	$\nu(V)=1$.  For $h>0$ and $f\in C_b(X)$, the semigroup identity yields
	\begin{align}\label{eq:occupation-invariance}
		\int_VP_h^\eta f\dd\nu_T-\int_Xf\dd\nu_T
		&=\frac1T\left(\int_T^{T+h}P_t^\eta f(v)\dd t
		-\int_0^hP_t^\eta f(v)\dd t\right),\notag\\
		\left|\int_VP_h^\eta f\dd\nu_T-\int_Xf\dd\nu_T\right|
		&\le\frac{2h\|f\|_\infty}{T}.
	\end{align}
	Use \eqref{eq:restricted-kernel-passage} for the first integral and weak
	convergence for the second.  It follows that
	$\int_VP_h^\eta f\dd\nu=\int_Xf\dd\nu$, hence $\nu\in\mathcal J_\eta$.
	
	\emph{Step 4: a bound for every invariant law.}
	Let $\nu\in\mathcal J_\eta$. We derive its moment bound by truncation.
	Choose a smooth nonincreasing $\chi:[0,\infty)\to[0,1]$, with
	$\chi=1$ on $[0,1]$ and $\chi=0$ on $[2,\infty)$, and put
	\[
	\Theta_N(z)=\int_0^z\chi(y/N)\dd y,
	\qquad \Theta_N'(z)=\chi(z/N),\qquad
	\Theta_N''(z)=N^{-1}\chi'(z/N)\le0.
	\]
	For $x=\|v\|_V^2$, the quadratic-variation rate of $\PhiV(v_\eta)$ is
	\[
	\mathcal Q_\eta(v)=4p^2(1+x)^{2p-2}
	[a_\eta(v)^2x^2+\|G_V^*v\|_{\mathfrak U}^2]\ge0.
	\]
	Consequently,
	\begin{align}\label{eq:stationary-truncation-generator}
		\mathscr L_\eta(\Theta_N\circ\PhiV)
		&=\chi(\PhiV/N)\mathscr L_\eta\PhiV
		+\frac{\chi'(\PhiV/N)}{2N}\mathcal Q_\eta\notag\\
		&\le\chi(\PhiV/N)[C-c(1+\eta\rho)\PhiV].
	\end{align}
	Start the process from $\nu$.  The function $\Theta_N\circ\PhiV$ is bounded;
	its nonzero drift and martingale coefficients are supported where
	$\PhiV\le2N$, a bounded $V$-ball.  Thus its martingale has mean zero, and
	stationarity gives
	\begin{align*}
		0&=\int_V\mathscr L_\eta(\Theta_N\circ\PhiV)\dd\nu,\\
		c\int_V\chi(\PhiV/N)(1+\eta\rho)\PhiV\dd\nu
		&\le C\int_V\chi(\PhiV/N)\dd\nu\le C.
	\end{align*}
	Since $\chi(\PhiV/N)\uparrow1$, monotone convergence proves
	\begin{equation}\label{eq:all-moment}
		\int_V(1+\eta\rho(v))\PhiV(v)\,\nu(\dd v)\le C/c.
	\end{equation}
	\emph{Step 5: return to the original variables.}
	Pathwise uniqueness and homogeneity give the dilation identity
	\begin{equation}\label{eq:semigroup-dilation}
		\widehat P_t^\eta(f\circ D_\eta)(u)
		=P_t^\eta f(D_\eta u).
	\end{equation}
	It follows that $\mu\in\mathcal I_\eta$ if and only if
	$\nu=(D_\eta)_\#\mu\in\mathcal J_\eta$.  Finally,
	\begin{align*}
		&\int_V(1+\rho(u))(\eta^2+\|u\|_V^2)^p\,\mu(\dd u)\\
		&\qquad=\eta^{2p}\int_V(1+\eta\rho(v))(1+\|v\|_V^2)^p\,\nu(\dd v)
		\le\frac Cc\eta^{2p}.
	\end{align*}
	This proves existence in $\mathcal I_\eta$ and both stationary estimates.
\end{proof}
	
	\subsection{Compactness and concentration}
	The stationary moment estimate also yields compactness of the invariant-measure sets and concentration of the original stationary states near
	the origin.
	
	\begin{proposition}[Compact stationary sets and concentration]\label{prop:stationary-set}
		Under the hypotheses of Theorem~\ref{thm:invariant}, for each fixed
		$\eta\in(0,1]$, both $\cJ_\eta$ and $\cI_\eta$ are nonempty, convex, and
		compact for weak convergence on $X$. The union
		$\bigcup_{0<\eta\le1}\cJ_\eta$ is relatively compact. Moreover, for every
		$\varepsilon>0$,
		\begin{equation}\label{eq:unscaled-concentration}
			\sup_{\mu\in\cI_\eta}\mu\{\norm{u}_V>\varepsilon\}
			\le \frac{C_{\mathrm{st}}\eta^q}{(\eta^2+\varepsilon^2)^p},\qquad
			\sup_{\mu\in\cI_\eta}\int\norm{u}_V^q\dd\mu
			\le C_{\mathrm{st}}\eta^q.
		\end{equation}
		Consequently, for every $f\in C_b(V)$,
		\[
		\lim_{\eta\downarrow0}\sup_{\mu\in\mathcal I_\eta}
		\left|\int_V f(u)\,\mu(\dd u)-f(0)\right|=0.
		\]
	\end{proposition}
	\begin{proof}
		The estimate \eqref{eq:all-moment}, together with compactness of $K_R$ in $X$,
		gives uniform tightness of $\bigcup_{0<\eta\le1}\mathcal J_\eta$. For fixed
		$\eta>0$, if
		$\nu_n\in\mathcal J_\eta$ and $\nu_n\Rightarrow\nu$, then
		\[
		\nu(V)=1,\qquad \int\PhiV\dd\nu\le C_{\mathrm{st}},\qquad
		\int P_t^\eta f\dd\nu
		=\lim_n\int P_t^\eta f\dd\nu_n
		=\lim_n\int f\dd\nu_n=\int f\dd\nu.
		\]
		The first equality of integrals follows from
		\eqref{eq:restricted-kernel-passage}.  Thus $\mathcal J_\eta$ is closed and
		compact.  Its convexity follows from
		$(a\nu+(1-a)\widetilde\nu)P_t^\eta=a\nu+(1-a)\widetilde\nu$ for
		$0\le a\le1$; dilation gives the same properties for $\mathcal I_\eta$.
		Moreover,
		\[
		(\eta^2+\varepsilon^2)^p\mu\{\|u\|_V>\varepsilon\}
		\le\int(\eta^2+\|u\|_V^2)^p\dd\mu
		\le C_{\mathrm{st}}\eta^q,
		\qquad
		\int\|u\|_V^q\dd\mu\le C_{\mathrm{st}}\eta^q.
		\]
		For $f\in C_b(V)$, split the integral over $\|u\|_V\le\varepsilon$ and its
		complement:
		\[
		\sup_{\mu\in\mathcal I_\eta}\left|\int f\dd\mu-f(0)\right|
		\le\sup_{\|u\|_V\le\varepsilon}|f(u)-f(0)|
		+\frac{2\|f\|_\infty C_{\mathrm{st}}\eta^q}
		{(\eta^2+\varepsilon^2)^p}.
		\]
		Let $\eta\downarrow0$ and then $\varepsilon\downarrow0$.
	\end{proof}

	\medskip

	\section{Stationary Cauchy-law limit}\label{sec:scaling}
	We identify the limit of the rescaled invariant family $\mathcal J_\eta$. First
	we compute the invariant law $\nu_*$ of the limiting linear equation. We then
	compare the nonlinear stationary family with $\nu_*$ by a localized common-noise
	coupling.
	
	For probability measures on $X$, let
	\begin{equation}\label{eq:BL}
		d_{\mathrm{BL},X}(\nu_1,\nu_2)
		=\sup_{\|f\|_\infty\le1,\ \Lip_X(f)\le1}
		\left|\int f\dd\nu_1-\int f\dd\nu_2\right|.
	\end{equation}
	
	\subsection{The limiting linear equation}
	We first identify the invariant law of the limiting linear equation
	obtained in the small-forcing limit.
	
	\begin{theorem}[Invariant law of the limiting linear equation]\label{thm:linear-limit}
		The linear equation defines Markov semigroups on both $V$ and $X$. In either
		state space, its unique invariant probability is supported in $V$ and is given by
		\begin{equation}\label{eq:ratio}
			\nu_* = \Law\left(\frac{Z}{a_0|\zeta|}\right),\qquad
			Z\sim\mathcal N(0,Q_V),\quad \zeta\sim\mathcal N(0,1),
		\end{equation}
		where $Z$ is a $V$-valued Gaussian and $Z,\zeta$ are independent. Under the
		embedding $V\hookrightarrow X$, the covariance of $Z$ is $Q_X$. For
		$\ell\in X^*$,
		\begin{equation}\label{eq:char}
			\int_Xe^{it\ell(v)}\nu_*(\dd v)
			=\exp\left(-\frac{|t|}{a_0}\norm{G_X^*\ell}_{\mathfrak U}\right).
		\end{equation}
		For every $0<\alpha<1$,
		\begin{equation}\label{eq:frac-moment}
			\int\norm{v}_X^\alpha\nu_*(\dd v)
			=\frac{2^{-\alpha/2}\Gamma((1-\alpha)/2)}{\sqrt\pi\,a_0^\alpha}
			\E\norm{Z}_X^\alpha<\infty.
		\end{equation}
		If $G\ne0$, every norm moment of order at least one is infinite; if $G=0$,
		then $\nu_*=\delta_0$.
	\end{theorem}
	\begin{proof}
		For $v\in V$, the finite-time solution has the adapted representation
		\begin{equation}\label{eq:adapted-linear}
			V_t^v=\mathcal E_t\left(v+\int_0^t\mathcal E_s^{-1}G\dd B_s\right),
			\qquad \mathcal E_t=\exp(a_0W_t-a_0^2t/2).
		\end{equation}
		The stochastic integral is defined by localization on each finite interval.
		Independence of $B$ and $W$ gives zero cross-variation.
		For two-sided noises, the stationary value is represented by
		\begin{equation}\label{eq:stationary-convolution}
			V_t^*=\int_{-\infty}^t
			e^{a_0(W_t-W_s)-a_0^2(t-s)/2}G\dd B_s.
		\end{equation}
		Conditioned on the full path of $W$, the integral is a Wiener integral with
		respect to the independent process $B$; equivalently, it is the limit of adapted
		solutions started at finite negative times.
		
		The conditional covariance is $I Q_V$, and its trace is $q_V I$, with
		\[
		I=\int_0^\infty e^{2a_0\widetilde W_r-a_0^2r}\dd r<\infty
		\quad\hbox{almost surely},
		\]
		where $\widetilde W$ is the Brownian motion obtained from the backward
		increments of $W$. Brownian scaling and Dufresne's identity, in the normalization
		of \cite{MY}, give
		\begin{equation}\label{eq:dufresne}
			I\ \stackrel{d}{=}\ \frac1{2a_0^2\Gamma_{1/2}}
			\ \stackrel{d}{=}\ \frac1{a_0^2\zeta^2}.
		\end{equation}
		Thus the conditional Gaussian mixture has the law \eqref{eq:ratio}. Shift
		invariance of the two-sided noise increments gives stationarity, and restarting
		the finite-past equation identifies it as a solution of
		\eqref{eq:linear-limit}.
		
		For common driving noises,
		\begin{equation}\label{eq:linear-limit-diff}
			V_t^v-V_t^w=\mathcal E_t(v-w),\qquad \mathcal E_t\longrightarrow0
			\quad\hbox{almost surely}.
		\end{equation}
		Couple any initial law with $\nu_*$ independently of the future noises. For
		bounded Lipschitz test functions, bounded convergence shows that the
		bounded--Lipschitz distance between the two time-$t$ laws tends to zero. Applying
		this to an invariant initial law proves uniqueness among all invariant
		probabilities. Formula~\eqref{eq:char} follows because a centered real Gaussian divided by an
		independent $|\zeta|$ is Cauchy, and \eqref{eq:frac-moment} follows by
		independence and a Gaussian integral. When $G\ne0$, some continuous linear
		functional has a nondegenerate Cauchy marginal; hence every norm moment of order
		at least one diverges.
	\end{proof}
	
	\subsection{Convergence of invariant measures}
	For $0<\alpha<1$, let
	\begin{equation}\label{eq:Walpha}
		\mathsf W_\alpha(\nu_1,\nu_2)
		=\inf_{\pi\in\Pi(\nu_1,\nu_2)}
		\int_{X\times X}\|x-y\|_X^\alpha\,\pi(\dd x,\dd y),
	\end{equation}
	where $\Pi(\nu_1,\nu_2)$ is the set of couplings.  Since
	$d_\alpha(x,y)=\|x-y\|_X^\alpha$ is a metric for $0<\alpha<1$,
	$\mathsf W_\alpha$ is its $1$-Wasserstein distance.
	
	The next theorem is the principal small-forcing result: it gives one common limit for the
	entire rescaled invariant-measure family, without assuming uniqueness at fixed positive forcing.
	
	\begin{theorem}[Uniform convergence of rescaled invariant measures]\label{thm:scaling}
		Assume the hypotheses of Theorem~\ref{thm:invariant}.  Let $G\ne0$ and $q=2p$.  Put
		\begin{equation}\label{eq:rate-exponent}
			\lambda_q=\frac{q(1-q)a_0^2}{2},\qquad
			\beta=1+2a_0^2,\qquad
			\vartheta=\left(1+\frac2q+\frac{\beta}{2\lambda_q}\right)^{-1}.
		\end{equation}
		There is a constant $K$, independent of the invariant measure, such that, for
		$0<\eta\le e^{-1}$,
		\begin{equation}\label{eq:rate}
			\sup_{\nu\in\mathcal J_\eta}
			d_{\mathrm{BL},X}(\nu,\nu_*)
			\le K(1+\log\eta^{-1})\eta^\vartheta.
		\end{equation}
		For every selection $\nu_\eta\in\mathcal J_\eta$, the corresponding stationary
		path laws converge in $C([0,T];X)$ to the stationary path law of
		\eqref{eq:linear-limit} on each finite time interval as $\eta\downarrow0$.
		Moreover, for every $0<\alpha<q$,
		\begin{equation}\label{eq:uniform-lower-moments}
			\sup_{\nu\in\mathcal J_\eta}
			\left|\int_X\|x\|_X^\alpha\,\nu(\dd x)
			-\int_X\|x\|_X^\alpha\,\nu_*(\dd x)\right|
			+\sup_{\nu\in\mathcal J_\eta}\mathsf W_\alpha(\nu,\nu_*)
			\longrightarrow0
			\qquad\text{as }\eta\downarrow0.
		\end{equation}
	\end{theorem}
	\begin{proof}
		\emph{Step 1: fractional contraction of the linear equation.}
		For $0<\alpha<1$, put
		$\lambda_\alpha=\alpha(1-\alpha)a_0^2/2$.  The common-noise identity
		\eqref{eq:linear-limit-diff}, together with
		$\E\mathcal E_t^\alpha=e^{-\lambda_\alpha t}$, gives, for every law $\nu$
		with finite $\alpha$th $X$-moment,
		\begin{equation}\label{eq:mix}
			d_{\mathrm{BL},X}((P_t^0)^*\nu,\nu_*)
			\le2^{1-\alpha}e^{-\lambda_\alpha t}
			\left(\int\norm{v}_X^\alpha\dd\nu+
			\int\norm{v}_X^\alpha\dd\nu_*\right).
		\end{equation}
		Here we used $2\wedge r\le2^{1-\alpha}r^\alpha$ for $r\ge0$ and
		$\norm{x-y}_X^\alpha\le\norm{x}_X^\alpha+\norm{y}_X^\alpha$.
		
		\emph{Step 2: localized common-noise comparison.}
		Start $v_\eta$ and $V$ from the same almost surely $V$-valued initial variable
		and drive them by the same pair $(W,B)$. Define
		\[
		\tau_R=\inf\{t\ge0:\|v_\eta(t)\|_V\ge R\},
		\]
		with $\tau_R=0$ when the initial state lies outside $K_R$. For $R\ge1$ and
		$\eta R\le1$,
		\begin{align}
			\E\norm{(v_\eta-V)(T\wedge\tau_R)}_X^2
			&\le C\eta^2R^4e^{\beta T},\qquad \beta=1+2a_0^2,
			\label{eq:point-compare}\\
			\E\sup_{t\le T}\norm{(v_\eta-V)(t\wedge\tau_R)}_X^2
			&\le C_T\eta^2R^4.\label{eq:path-compare}
		\end{align}
		For $w=v_\eta-V$ before $\tau_R$,
		\[
		\dd w=-\eta F(v_\eta)\dd t+[a_0w+g_\eta(v_\eta)]\dd W_t,
		\qquad g_\eta(v)=[\psi(\eta\rho(v))-a_0]v.
		\]
		The additive noises cancel. On the stopped ball,
		\begin{equation}\label{eq:forcing-error}
			\eta\norm{F(v)}_X+\norm{g_\eta(v)}_X\le C\eta R^2.
		\end{equation}
		Before $\tau_R$, It\^o's formula gives
		\begin{align*}
			\dd\|w\|_X^2
			={}&[-2\eta\ip{F(v_\eta)}{w}_X
			+\|a_0w+g_\eta(v_\eta)\|_X^2]\dd t
			+2\ip{w}{a_0w+g_\eta(v_\eta)}_X\dd W_t\\
			\le{}&[\beta\|w\|_X^2+C\eta^2R^4]\dd t+\dd M_t,\\
			\dd[M]_t\le{}&C(\|w\|_X^4+\eta^2R^4\|w\|_X^2)\dd t.
		\end{align*}
		Since $w(0)=0$, localization and Gronwall give
		\[
		\E\|w(T\wedge\tau_R)\|_X^2
		\le C\eta^2R^4\int_0^T e^{\beta(T-s)}\dd s
		\le C\eta^2R^4e^{\beta T}.
		\]
		The BDG estimate for the displayed bracket, followed by Young's inequality,
		absorbs one half of $\E\sup_{t\le T}\|w(t\wedge\tau_R)\|_X^2$ and yields
		\eqref{eq:path-compare}. The bounds hold conditionally on the common initial
		state, with constants uniform on $K_R$. By the convention $\tau_R=0$ outside
		$K_R$, the stopped difference vanishes there, so the estimates integrate against
		any initial law.
		
		\emph{Step 3: marginal laws and the rate.}
		By Theorem~\ref{thm:invariant},
		$\int\PhiV\dd\nu_\eta\le C_{\mathrm{st}}$ uniformly in
		$\nu_\eta\in\mathcal J_\eta$. Start $v_\eta$ and $V$ from a common initial
		variable with law $\nu_\eta$. Equation~\eqref{eq:exit} gives
		\[
		\PP(\tau_R\le T)\le (C_{\mathrm{st}}+CT)(1+R^2)^{-p}.
		\]
		On $\{\tau_R>T\}$ use \eqref{eq:point-compare}; on its complement, use the
		bound $2$ for the difference of two bounded test functions. For every $f$ in the
		defining class of \eqref{eq:BL},
		\[
		\left|\int(P_T^\eta f-P_T^0f)\dd\nu_\eta\right|
		\le C\eta R^2e^{\beta T/2}+C(1+T)R^{-q}.
		\]
		Invariance of $\nu_\eta$ and the linear contraction \eqref{eq:mix} therefore imply the three-error estimate
		\begin{equation}\label{eq:three-errors}
			d_{\mathrm{BL},X}(\nu_\eta,\nu_*)
			\le K\left[\eta R^2e^{\beta T/2}+(1+T)R^{-q}
			+e^{-\lambda_qT}\right].
		\end{equation}
		Choose
		\[
		R=\eta^{-\vartheta/q},\qquad T=\frac{\vartheta}{\lambda_q}\log\eta^{-1}.
		\]
		Since $\vartheta/q<1$, the constraint $\eta R\le1$ holds. The identity
		$1-2\vartheta/q-\beta\vartheta/(2\lambda_q)=\vartheta$ then gives
		\eqref{eq:rate}.
		
		\emph{Step 4: stationary paths.}
		For fixed $T$, the path comparison and the exit bound give
		\begin{equation}\label{eq:path-bounded}
			\E\left[1\wedge\sup_{t\le T}\norm{v_\eta(t)-V(t)}_X^2\right]
			\le C_T\eta^2R^4+(C_{\mathrm{st}}+CT)(1+R^2)^{-p}.
		\end{equation}
		Send $\eta\downarrow0$ with $R$ fixed, and then $R\to\infty$.
		The linear path kernel is continuous from $X$ to probability measures on
		$C([0,T];X)$: by \eqref{eq:linear-limit-diff},
		\[
		\sup_{t\le T}\norm{V_t^{x_n}-V_t^x}_X
		\le\sup_{t\le T}\mathcal E_t\norm{x_n-x}_X\longrightarrow0
		\quad\hbox{almost surely}
		\]
		whenever $x_n\to x$. Weak convergence
		$\nu_\eta\Rightarrow\nu_*$ and continuity of this kernel therefore imply
		convergence of the corresponding linear path laws. Combining this with
		\eqref{eq:path-bounded} proves path-law convergence.
		
		\emph{Step 5: lower moments.}
		For $0<\alpha<q$, the continuous embedding $V\hookrightarrow X$ gives
		\begin{equation}\label{eq:lower-moment-tail}
			\sup_{\substack{0<\eta\le1\\\nu\in\mathcal J_\eta}}
			\int_{\{\|x\|_X>R\}}\|x\|_X^\alpha\,\nu(\dd x)
			\le R^{\alpha-q}
			\sup_{\eta,\nu}\int\|x\|_X^q\,\nu(\dd x)
			\le C R^{\alpha-q}.
		\end{equation}
		The same bound holds for $\nu_*$.  Hence weak convergence applied to
		$\|x\|_X^\alpha\wedge R^\alpha$, followed by $R\to\infty$, gives moment
		convergence.  For any sequence $\eta_n\downarrow0$ and
		$\nu_n\in\mathcal J_{\eta_n}$, a Skorokhod coupling gives
		$Z_n\to Z_*$ almost surely, with laws $\nu_n$ and $\nu_*$.  Since
		\[
		\sup_n\E\|Z_n-Z_*\|_X^q
		\le\sup_n\E\|Z_n\|_X^q+\E\|Z_*\|_X^q<\infty,
		\]
		the variables $\|Z_n-Z_*\|_X^\alpha$ are uniformly integrable.  Therefore
		\[
		\mathsf W_\alpha(\nu_n,\nu_*)
		\le\E\|Z_n-Z_*\|_X^\alpha\longrightarrow0.
		\]
		Because the argument applies to every sequence $\eta_n\downarrow0$, both
		suprema in \eqref{eq:uniform-lower-moments} tend to zero.
	\end{proof}
	
	The uniform convergence implies both collapse of the stationary set and convergence of its
	low-order moments.
	
	\begin{corollary}[Vanishing diameter and moment scaling]\label{cor:vanishing-diameter}
		Under Theorem~\ref{thm:scaling}, with $G\ne0$, the scaled invariant family
		satisfies
		\begin{equation}\label{eq:diameter}
			\sup_{\nu,\widetilde\nu\in\cJ_\eta}
			d_{\mathrm{BL},X}(\nu,\widetilde\nu)
			\le2K(1+\log\eta^{-1})\eta^\vartheta.
		\end{equation}
		For $0<\alpha<q$, set
		\[
		m_\alpha=\int_X\|v\|_X^\alpha\,\nu_*(\dd v)
		=\frac{2^{-\alpha/2}\Gamma((1-\alpha)/2)}{\sqrt\pi\,a_0^\alpha}
		\E\|Z\|_X^\alpha\in(0,\infty).
		\]
		Then
		\begin{equation}\label{eq:exact-stationary-amplitude}
			\lim_{\eta\downarrow0}\sup_{\mu\in\cI_\eta}
			\left|\eta^{-\alpha}\int_X\|u\|_X^\alpha\,\mu(\dd u)-m_\alpha\right|=0.
		\end{equation}
		In particular, for sufficiently small $\eta>0$,
		\begin{equation}\label{eq:two-sided-stationary-amplitude}
			\frac{m_\alpha}{2}\eta^\alpha
			\le\inf_{\mu\in\cI_\eta}\int_X\|u\|_X^\alpha\,\mu(\dd u)
			\le\sup_{\mu\in\cI_\eta}\int_X\|u\|_X^\alpha\,\mu(\dd u)
			\le\frac{3m_\alpha}{2}\eta^\alpha.
		\end{equation}
		Stationary path convergence on each fixed finite time window is uniform over all
		choices of invariant initial law, in a bounded--Lipschitz metric for the uniform
		$X$-path topology.
	\end{corollary}
	\begin{proof}
		The diameter bound is the triangle inequality. Theorem~\ref{thm:scaling} gives
		uniform convergence of the $\alpha$th $X$-moments, and dilation gives
		\[
		\int_X\|v\|_X^\alpha\,(D_\eta)_\#\mu(\dd v)
		=\eta^{-\alpha}\int_X\|u\|_X^\alpha\,\mu(\dd u).
		\]
		This proves \eqref{eq:exact-stationary-amplitude}. The formula for $m_\alpha$
		follows from \eqref{eq:frac-moment}. Since $G_X\ne0$, a continuous functional
		detects a nonzero Gaussian direction, so $m_\alpha>0$. Taking the uniform error
		smaller than $m_\alpha/2$ yields \eqref{eq:two-sided-stationary-amplitude}.
		If uniform path convergence failed, one could choose $\eta_n\downarrow0$ and a
		sequence of invariant laws for which convergence fails. The path comparison
		\eqref{eq:path-bounded} and continuity of the linear path kernel apply with the
		same constants, giving a contradiction.
	\end{proof}
	
	\begin{remark}[Topology, moments, and rate]
		The moment and Wasserstein conclusions hold in $X$ for orders $\alpha<q$.
		This restriction is consistent with the Cauchy tails of $\nu_*$, whose
		nondegenerate one-dimensional marginals have infinite moments of order at least
		one. The exponent $\vartheta$ is obtained by
		balancing the comparison, exit, and linear-contraction errors in
		\eqref{eq:three-errors}.
	\end{remark}
	
	\begin{remark}[Beyond exact quadratic homogeneity]\label{rem:asymptotic}
		If the local theory, the high-order bound \eqref{eq:b}, and the resulting
		scale-uniform stationary and exit estimates remain valid, exact quadratic
		homogeneity can be relaxed in the comparison step.  It suffices that, for every
		fixed $R$,
		\begin{equation}\label{eq:asymptotic}
			\varepsilon_\eta(R)=
			\sup_{v\in K_R}\|\eta^{-1}F(\eta v)\|_X\longrightarrow0.
		\end{equation}
		Then the comparison error $C\eta R^2$ is replaced by
		\[
		\varepsilon_\eta(R)+R\sup_{0\le z\le C_\rho\eta R}|\psi(z)-a_0|.
		\]
		This gives qualitative convergence by first choosing $T$, then $R$, and finally
		$\eta$.  An algebraic rate additionally requires uniform bounds on radii
		$R=R(\eta)$.
	\end{remark}

	\medskip 
	
	\section{Hamiltonian transport models}\label{sec:models}
	For each model we specify the phase spaces and feedback functional, prove the
	lower-space difference estimate, identify the high-order energy pairing, and verify
	the regularization conditions of Theorem~\ref{thm:local-criterion}. The abstract
	global and stationary results then follow from
	Theorems~\ref{thm:global},~\ref{thm:invariant}, and~\ref{thm:scaling}.
	The scaling limit is nontrivial
	when $G\ne0$; for $G=0$, Theorem~\ref{thm:global} gives the invariant law
	$\delta_0$. We first record three Sobolev estimates used throughout.
	Appendix~\ref{sec:low} uses the conserved low energies identified for each
	model.
	
	\begin{lemma}[Product, commutator, and transport estimates]\label{lem:sobolev}
		Let $s>0$, $\Lambda^s=(I-\Delta)^{s/2}$, and let the functions below be
		smooth on $\T^d$.  The estimates extend to Sobolev functions for which their
		right-hand sides are finite.  For scalar functions $f,g$,
		\begin{equation}\label{eq:sobolev-product}
			\|fg\|_{H^s}\le C_s
			(\|f\|_{L^\infty}\|g\|_{H^s}
			+\|g\|_{L^\infty}\|f\|_{H^s}).
		\end{equation}
		For a vector field $f$ and a scalar or vector field $w$,
		\begin{equation}\label{eq:KP}
			\|[\Lambda^s,f\cdot\nabla]w\|_{L^2}
			\le C_s\big(\|\nabla f\|_{L^\infty}\|w\|_{H^s}
			+\|f\|_{H^s}\|\nabla w\|_{L^\infty}\big),
		\end{equation}
		and
		\begin{align}\label{eq:transport-energy}
			\big|\ip{\Lambda^s(f\cdot\nabla w)}{\Lambda^sw}_{L^2}\big|
			\le C_s\big(&\|\nabla f\|_{L^\infty}\|w\|_{H^s}^2\notag\\
			&+\|f\|_{H^s}\|\nabla w\|_{L^\infty}\|w\|_{H^s}\big).
		\end{align}
	\end{lemma}
	\begin{proof}
		The periodic product and Kato--Ponce estimates give
		\eqref{eq:sobolev-product}--\eqref{eq:KP}; see \cite{Li}.  Integration by
		parts gives the identity
		\[
		\ip{\Lambda^s(f\cdot\nabla w)}{\Lambda^sw}_{L^2}
		=-\frac12\int_{\T^d}(\nabla\cdot f)|\Lambda^sw|^2\dd x
		+\ip{[\Lambda^s,f\cdot\nabla]w}{\Lambda^sw}_{L^2}.
		\]
		Cauchy--Schwarz and \eqref{eq:KP} prove \eqref{eq:transport-energy}.
	\end{proof}
	
	\begin{lemma}[Bilinear regularization residual]\label{lem:bilinear-residual}
		Let $E^a$ denote scalar, vector, or finite-product periodic Sobolev spaces.
		Fix $r>d/2$ and $s>r+1$. Let
		$\mathcal B:E^r\times E^{r+1}\to E^r$ be a bilinear map and set
		$F(v)=\mathcal B(v,v)$ for $v\in E^s$. Assume that
		\begin{equation}\label{eq:bilinear-bounds}
			\|\mathcal B(f,g)\|_{E^r}\le C\|f\|_{E^r}\|g\|_{E^{r+1}},
			\qquad \|F(v)\|_{E^{s-1}}\le C\|v\|_{E^s}^2.
		\end{equation}
		For $F_n(v)=J_nF(J_nv)$ and $\delta=s-r-1$, one has
		\begin{equation}\label{eq:bilinear-residual-bound}
			\sup_{\|v\|_{E^s}\le M}\|F_n(v)-F(v)\|_{E^r}
			\le C_Mn^{-\delta}.
		\end{equation}
		Moreover, $F$ is continuous in the $E^r$ topology on each bounded $E^s$-ball.
		The estimates may be applied in the ambient product space when the final drift
		and regularizers preserve a constraint subspace.
	\end{lemma}
	\begin{proof}
		Put $e_n=(J_n-I)v$. Bilinearity gives the full decomposition
		\begin{equation}\label{eq:three-residuals}
			F_n(v)-F(v)
			=J_n\mathcal B(e_n,J_nv)+J_n\mathcal B(v,e_n)
			+(J_n-I)F(v).
		\end{equation}
		The three terms satisfy, respectively,
		\begin{align*}
			\|J_n\mathcal B(e_n,J_nv)\|_{E^r}
			&\le C\|e_n\|_{E^r}\|v\|_{E^{r+1}}
			\le C_Mn^{-(s-r)},\\
			\|J_n\mathcal B(v,e_n)\|_{E^r}
			&\le C\|v\|_{E^r}\|e_n\|_{E^{r+1}}
			\le C_Mn^{-\delta},\\
			\|(J_n-I)F(v)\|_{E^r}
			&\le Cn^{-\delta}\|F(v)\|_{E^{s-1}}
			\le C_Mn^{-\delta}.
		\end{align*}
		For $v,w$ in that ball,
		\[
		\|F(v)-F(w)\|_{E^r}
		\le C_M\bigl(\|v-w\|_{E^r}+\|v-w\|_{E^{r+1}}\bigr).
		\]
		Interpolation with $\theta=(s-r-1)/(s-r)\in(0,1)$ gives
		$\|v-w\|_{E^{r+1}}\le(2M)^{1-\theta}\|v-w\|_{E^r}^{\theta}$,
		which proves the asserted continuity and verifies (R3).
	\end{proof}
	
	\begin{lemma}[Recovery of transport energy and paired cancellations]\label{lem:transport-recovery}
		Let $s>d/2+1$, $a\in H^s(\T^d;\R^d)$, and let $w,z$ be scalar or vector
		fields in $H^s$. Put $w_s=\Lambda^sw$, $z_s=\Lambda^sz$ and
		$C_s(a,w)=[\Lambda^s,a\cdot\nabla]w$. Then
		\begin{align}\label{eq:transport-recovery-identity}
			\ip{J_n\Lambda^s(a\cdot\nabla w)}{J_nw_s}_{L^2}
			={}&-\frac12\int(\nabla\cdot a)|J_nw_s|^2\dd x\notag\\
			&+\ip{[J_n,a\cdot\nabla]w_s}{J_nw_s}_{L^2}
			+\ip{J_nC_s(a,w)}{J_nw_s}_{L^2}.
		\end{align}
		The left side converges to
		$-\frac12\int(\nabla\cdot a)|w_s|^2+\ip{C_s(a,w)}{w_s}_{L^2}$,
		and is bounded uniformly in $n$ by the right side of
		\eqref{eq:transport-energy}. If $\nabla\cdot a=0$, then the two cross
		principal terms cancel exactly at every regularization level:
		\begin{equation}\label{eq:regularized-paired-cancellation}
			\ip{a\cdot\nabla J_nw_s}{J_nz_s}_{L^2}
			+\ip{a\cdot\nabla J_nz_s}{J_nw_s}_{L^2}=0.
		\end{equation}
		The remaining terms in the paired energy converge separately by the same
		commutator argument.
	\end{lemma}
	\begin{proof}
		In distributions,
		\[
		J_n\Lambda^s(a\cdot\nabla w)
		=a\cdot\nabla J_nw_s+[J_n,a\cdot\nabla]w_s+J_nC_s(a,w).
		\]
		Pairing with $J_nw_s$ and integrating the principal term by parts gives
		\eqref{eq:transport-recovery-identity}. The middle term tends to zero by
		\eqref{eq:mollifier-comm}. The last term converges since
		$C_s(a,w)\in L^2$ by \eqref{eq:KP}. The same bounds and the contraction
		property of $J_n$ give the uniform energy estimate. Pairing instead with
		$J_nz_s$, adding the expression with $w,z$ exchanged, and integrating by parts
		gives \eqref{eq:regularized-paired-cancellation}. The relevant pairings are
		continuous functions of $(a,w,z)$ in $H^s$, by the Sobolev embedding and the
		bilinear commutator bound. Thus the energy forms extend continuously to the
		stated Sobolev phase space.
	\end{proof}
	
	\subsection{Active scalars}\label{sec:active}
	On the mean-zero subspace of $\T^2$, let $\Lambda_0=(-\Delta)^{1/2}$:
	for every real $\gamma$,
	\[
	\widehat{\Lambda_0^\gamma f}(k)=|k|^\gamma\widehat f(k)
	\quad(k\ne0),\qquad
	\widehat{\Lambda_0^\gamma f}(0)=0.
	\]
	This homogeneous diagnostic multiplier is distinct from
	$\Lambda=(I-\Delta)^{1/2}$ in the Sobolev norm. Define
	\begin{equation}\label{eq:Kalpha}
		\mathcal K_\alpha=\nabla^\perp\Lambda_0^{-\alpha},\qquad
		1\le\alpha\le2,
		\qquad F_\alpha(\theta)=(\mathcal K_\alpha\theta)\cdot\nabla\theta.
	\end{equation}
	We fix the mean-zero scalar space and set
	\begin{equation}\label{eq:active-space}
		V=H^s_0(\T^2),\quad X=H^r_0(\T^2),\quad
		s>3,\quad 2<r<s-1.
	\end{equation}
	Use the feedback
	\begin{equation}\label{eq:active-rho}
		\rho_\alpha(\theta)=\norm{\theta}_{W^{1,\infty}}
		+\norm{\mathcal K_\alpha\theta}_{W^{1,\infty}}.
	\end{equation}
	Up to the convention for $\nabla^\perp$, the endpoints give SQG and
	two-dimensional Euler vorticity dynamics. For the SQG model and more singular
	active scalars, see \cite{CMT,CCCGW}. Transport-noise Euler solutions and
	stationary enstrophy scaling are treated in \cite{LC,FL20}; deterministic Euler
	invariant measures in \cite{Latocca23,BFerrari}; and modified SQG invariant
	measures in \cite{NPST18}. Here $1\le\alpha\le2$, so the diagnostic velocity
	operator has nonpositive order.
	
	We begin the model verification with active scalars, whose nonlinear drift is a scalar transport
	term. 
	\begin{proposition}[Active-scalar verification]\label{prop:active}
		For \eqref{eq:Kalpha}--\eqref{eq:active-rho}, the drift and regularizers
		satisfy Assumption~\ref{ass:local-structure} and the hypotheses of
		Theorem~\ref{thm:local-criterion}. Thus every locally Lipschitz feedback gives a
		unique maximal $V$-valued, $X$-stable local solution. If
		Assumption~\ref{ass:stabilization} holds with the corresponding energy constant,
		Theorems~\ref{thm:global},~\ref{thm:invariant}, and~\ref{thm:scaling}
		apply. The conserved low energy is
		$E_0(\theta)=\norm{\theta}_{L^2}^2$.  For square-root feedback and
		$G\ne0$, Theorem~\ref{thm:low} gives an infinite stationary $H^s$ second
		moment at every fixed $\eta>0$.
	\end{proposition}
	\begin{proof}
		For $\alpha\ge1$, $\mathcal K_\alpha$ has order $1-\alpha\le0$, is bounded on each $H^j_0$, commutes with Fourier regularizers, and is divergence free. Since $H^r(\T^2)\hookrightarrow W^{1,\infty}$, \eqref{eq:rho} holds for \eqref{eq:active-rho}. The commutator estimate gives
		\begin{align*}
			\big|\ip{\Lambda^sF_\alpha(\theta)}{\Lambda^s\theta}\big|
			&\le C_s\big(\norm{\nabla\mathcal K_\alpha\theta}_\infty\norm{\theta}_{H^s}
			+\norm{\mathcal K_\alpha\theta}_{H^s}\norm{\nabla\theta}_\infty\big)
			\norm{\theta}_{H^s}\\
			&\le C_{s,\alpha}\rho_\alpha(\theta)\norm{\theta}_{H^s}^2.
		\end{align*}
		The feedback \eqref{eq:active-rho} therefore controls both commutator terms.
		
		For $w=\theta_1-\theta_2$, decompose the drift difference as
		\[
		(\mathcal K_\alpha\theta_1)\cdot\nabla w
		+(\mathcal K_\alpha w)\cdot\nabla\theta_2.
		\]
		The first term is estimated in the $H^r$ energy using \eqref{eq:transport-energy} and $H^r\hookrightarrow W^{1,\infty}$. The second satisfies
		\[
		\norm{(\mathcal K_\alpha w)\cdot\nabla\theta_2}_{H^r}
		\le C\norm{w}_{H^r}\norm{\theta_2}_{H^{r+1}}
		\le C_R\norm{w}_{H^r}.
		\]
		This proves \eqref{eq:diff}, while the product estimate proves
		$\norm{F_\alpha(\theta)}_{H^r}\le C\norm{\theta}_{H^s}^2$.
		The drift is quadratic. Finally, incompressibility gives
		\[
		\ip{F_\alpha(\theta)}{\theta}_{L^2}
		=\frac12\int_{\T^2}\mathcal K_\alpha\theta\cdot\nabla(\theta^2)\dd x=0.
		\]
		Mean zero is preserved because the drift is a divergence and the noises take values in the mean-zero space.
		
		For the local construction, set
		$\mathcal B_\alpha(f,g)=(\mathcal K_\alpha f)\cdot\nabla g$ and
		$F_n(v)=J_nF_\alpha(J_nv)$. The Sobolev product estimate gives
		\[
		\|\mathcal B_\alpha(f,g)\|_{H^r}
		\le C\|f\|_{H^r}\|g\|_{H^{r+1}},\qquad
		\|F_\alpha(v)\|_{H^{s-1}}\le C\|v\|_{H^s}^2.
		\]
		Lemma~\ref{lem:bilinear-residual} therefore verifies (R2)--(R3), including the
		outer regularization term. For (R1), self-adjointness and commutation give
		\[
		2\ip{F_n(v)}{v}_{H^s}
		=2\ip{F_\alpha(J_nv)}{J_nv}_{H^s},\qquad
		\rho_\alpha(J_nv)\le C\rho_\alpha(v).
		\]
		Since $J_n$ has finite-dimensional range, $F_n:V\to V$ is locally Lipschitz.
		Finally, with $u=\mathcal K_\alpha\theta$ and $\theta_s=\Lambda^s\theta$,
		\[
		b_V(\theta)=-2\ip{C_s(u,\theta)}{\theta_s}_{L^2}.
		\]
		This form is continuous and cubic on $H^s_0$. Lemma~\ref{lem:transport-recovery},
		using $\nabla\cdot u=0$, gives the pointwise limit and the uniform domination
		in (R4). This verifies the local and energy hypotheses.
	\end{proof}
	
	The Lie--Poisson interpretation can be seen directly. Let
	$[f,g]=\partial_1f\partial_2g-\partial_2f\partial_1g$,
	$\mathcal J(\theta)f=[f,\theta]$, and
	$\mathcal H_\alpha(\theta)=\tfrac12\int\theta\Lambda_0^{-\alpha}\theta\dd x$.
	With the chosen sign convention,
	$F_\alpha=\mathcal J(\theta)D\mathcal H_\alpha(\theta)$. The $L^2$ energy
	used above is a Casimir of the same Lie--Poisson bracket, whereas
	$\mathcal H_\alpha$ is the Hamiltonian. The proof uses the nonpositive order of
	$\mathcal K_\alpha$ and therefore covers $1\le\alpha\le2$.
	
	\subsection{Ideal magnetohydrodynamics}\label{sec:mhd}
	Let $H^j_\sigma(\T^2;\R^2)$ denote the mean-zero divergence-free vector fields, and let $\mathbb P$ be the corresponding Leray projection. For $U=(u,b)$, define
	\begin{equation}\label{eq:mhd-drift}
		F_{\mathrm{MHD}}(u,b)=
		\begin{pmatrix}
			\mathbb P[(u\cdot\nabla)u-(b\cdot\nabla)b]\\
			\mathbb P[(u\cdot\nabla)b-(b\cdot\nabla)u]
		\end{pmatrix}.
	\end{equation}
	Here $u$ is the velocity and $b$ the magnetic field. Ideal incompressible MHD
	is a Lie--Poisson system on a semidirect product, with Hamiltonian equal to the
	total kinetic and magnetic energy \cite{HMR}; its transport form and
	Els\"asser variables are reviewed in \cite{CF}. Existing stochastic stationary
	results use viscous, tamed, or singularly perturbed models
	\cite{BDP07,FFGR17,Motyl22,Schenke22}, while transport noise provides a different
	regularization mechanism in three dimensions \cite{LuoMHD23}. We consider the
	inviscid two-dimensional system with common radial It\^o noise.
	
	Take
	\begin{equation}\label{eq:mhd-spaces}
		V=H^s_\sigma\times H^s_\sigma,\quad
		X=H^r_\sigma\times H^r_\sigma,\quad
		s>3,\quad 2<r<s-1,
	\end{equation}
	with the sum-of-squares norms, and set
	\begin{equation}\label{eq:mhd-rho}
		\rho_{\mathrm{MHD}}(U)=\norm{u}_{W^{1,\infty}}
		+\norm{b}_{W^{1,\infty}}.
	\end{equation}
	The scalar radial noise acts on the entire product state:
	\begin{equation}\label{eq:mhd-spde}
		\dd U_\eta+F_{\mathrm{MHD}}(U_\eta)\dd t
		=\psi(\rho_{\mathrm{MHD}}(U_\eta))U_\eta\dd W_t
		+\eta G\dd B_t,\qquad G\in\HS(\mathfrak U,V).
	\end{equation}
	The two components of $G$ may be correlated, and either component may be
	unforced. The same scalar Brownian motion $W$ multiplies both state components,
	which produces the radial energy bracket used in
	Section~\ref{sec:framework}.
	
	We next consider ideal MHD, where the key estimate is cancellation of the two principal
	cross-transport terms.
	
	\begin{proposition}[Ideal MHD]\label{prop:mhd}
		For \eqref{eq:mhd-spde} on the spaces \eqref{eq:mhd-spaces}, the drift and
		regularizers satisfy Assumption~\ref{ass:local-structure} and the hypotheses of
		Theorem~\ref{thm:local-criterion}. If Assumption~\ref{ass:stabilization} holds
		with the corresponding energy constant,
		Theorems~\ref{thm:global},~\ref{thm:invariant}, and~\ref{thm:scaling}
		apply. The positive conserved low energy is
		\begin{equation}\label{eq:mhd-energy}
			E_{\mathrm{MHD}}(u,b)=\norm{u}_{L^2}^2+\norm{b}_{L^2}^2.
		\end{equation}
		If $G\ne0$ and $\psi(y)=\sigma\sqrt{1+y}$ has the required strength, every
		invariant law at fixed $\eta>0$ has infinite product-$H^s$ second moment.
	\end{proposition}
	\begin{proof}
		\emph{Constraints and growth.}
		The Leray projection is self-adjoint and commutes with all scalar Fourier multipliers. The second component of \eqref{eq:mhd-drift} is divergence free: in its divergence the two contractions of $\nabla u$ and $\nabla b$ cancel after exchanging indices. Each transport term has zero spatial mean. The noises preserve both linear constraints because $G$ has range in $V$. Sobolev embedding makes \eqref{eq:mhd-rho} Lipschitz in $X$. Product estimates give the quadratic growth bound into $X$ and exact degree-two homogeneity.
		
		\emph{High-order energy.}
		Write $u_s=\Lambda^s u$ and $b_s=\Lambda^s b$. The principal diagonal terms vanish by $\nabla\cdot u=0$. The two magnetic principal terms cancel \emph{together}:
		\begin{equation}\label{eq:mhd-cross}
			-\ip{b\cdot\nabla b_s}{u_s}_{L^2}
			-\ip{b\cdot\nabla u_s}{b_s}_{L^2}
			=-\int_{\T^2}b\cdot\nabla(u_s\cdot b_s)\dd x=0.
		\end{equation}
		The four remaining terms are commutators
		$[\Lambda^s,u\cdot\nabla]u$, $[\Lambda^s,b\cdot\nabla]b$,
		$[\Lambda^s,u\cdot\nabla]b$, and
		$[\Lambda^s,b\cdot\nabla]u$, with the appropriate signs.
		By \eqref{eq:KP}, their total energy pairing is bounded by
		\[
		C_s\rho_{\mathrm{MHD}}(U)
		(\norm{u}_{H^s}^2+\norm{b}_{H^s}^2).
		\]
		This is the needed high-order form, with the factor two absorbed in $C_V$. Estimating the magnetic principal terms separately in $H^s$ would lose a derivative; \eqref{eq:mhd-cross} is essential.
		
		\emph{Lower-space stability.}
		Put $w=u_1-u_2$ and $z=b_1-b_2$. Since the Leray projection is orthogonal,
		we suppress it in the energy calculation. The drift difference equals
		\begin{align*}
			\delta F_u&=u_1\cdot\nabla w-b_1\cdot\nabla z
			+w\cdot\nabla u_2-z\cdot\nabla b_2,\\
			\delta F_b&=u_1\cdot\nabla z-b_1\cdot\nabla w
			+w\cdot\nabla b_2-z\cdot\nabla u_2.
		\end{align*}
		The first two terms in each line have the same diagonal and cross cancellation as above at level $r$. Their commutators are bounded by $C_M(\norm{w}_{H^r}^2+\norm{z}_{H^r}^2)$ on a product-$H^s$ ball. Each remaining term is controlled by, for example,
		\[
		\norm{w\cdot\nabla u_2}_{H^r}
		\le C\norm{w}_{H^r}\norm{u_2}_{H^{r+1}}.
		\]
		This gives \eqref{eq:diff} in the product space.
		
		\emph{Regularization residual and recovery of the coupled energy.}
		For $U_i=(u_i,b_i)$ define
		\[
		\mathcal B_M(U_1,U_2)=
		\begin{pmatrix}
			\mathbb P(u_1\cdot\nabla u_2-b_1\cdot\nabla b_2)\\
			\mathbb P(u_1\cdot\nabla b_2-b_1\cdot\nabla u_2)
		\end{pmatrix},\qquad F_{\mathrm{MHD}}(U)=\mathcal B_M(U,U).
		\]
		The componentwise product estimate yields
		\[
		\|\mathcal B_M(U_1,U_2)\|_{H^r\times H^r}
		\le C\|U_1\|_{H^r\times H^r}\|U_2\|_{H^{r+1}\times H^{r+1}},\qquad
		\|F_{\mathrm{MHD}}(U)\|_{H^{s-1}\times H^{s-1}}
		\le C\|U\|_V^2.
		\]
		Hence $F_n=J_nF_{\mathrm{MHD}}(J_n\cdot)$ satisfies (R2)--(R3) by
		Lemma~\ref{lem:bilinear-residual}. Self-adjointness, the high-order cancellation
		at the smooth state $J_nU$, and $\rho_{\mathrm{MHD}}(J_nU)\le C\rho_{\mathrm{MHD}}(U)$
		give (R1).
		
		To verify (R4) at a fixed $U\in V$, set
		\[
		C_{n,s}(a,w)=J_nC_s(a,w)+[J_n,a\cdot\nabla]\Lambda^sw.
		\]
		The diagonal principal terms vanish, and the magnetic principal terms satisfy
		\[
		-\ip{b\cdot\nabla J_nb_s}{J_nu_s}
		-\ip{b\cdot\nabla J_nu_s}{J_nb_s}=0.
		\]
		Thus, using also commutation of $\mathbb P$ and $J_n$,
		\begin{align}\label{eq:mhd-regularized-energy}
			\ip{J_nF_{\mathrm{MHD}}(U)}{J_nU}_V
			={}&\ip{C_{n,s}(u,u)-C_{n,s}(b,b)}{J_nu_s}_{L^2}\notag\\
			&+\ip{C_{n,s}(u,b)-C_{n,s}(b,u)}{J_nb_s}_{L^2}.
		\end{align}
		Lemma~\ref{lem:transport-recovery} gives
		$C_{n,s}(a,w)\to C_s(a,w)$ in $L^2$ and bounds the sum by
		$C_s\rho_{\mathrm{MHD}}(U)\|U\|_V^2$. Therefore (R4) holds with
		\begin{align*}
			b_V(U)=-2\big[&\ip{C_s(u,u)-C_s(b,b)}{u_s}_{L^2}
			+\ip{C_s(u,b)-C_s(b,u)}{b_s}_{L^2}\big].
		\end{align*}
		The derivative-losing cross terms therefore cancel \emph{before} the
		regularizer is removed. The displayed expression is continuous and cubic on
		$V$, which verifies the local hypotheses and the high-order energy identity.
		
		\emph{Low energy.}
		At level zero all commutators disappear, so
		$DE_{\mathrm{MHD}}(U)F_{\mathrm{MHD}}(U)=0$.
		The inclusion from $V$ into the product $L^2$ space is continuous and
		injective. Consequently, $G\ne0$ implies $q_0>0$, and square-root feedback
		satisfies \eqref{eq:upper-comparison}. The final assertion follows from
		Theorem~\ref{thm:low}.
	\end{proof}
	
	\subsection{Smoothing-inertia EP equations}\label{sec:epdiff}
	Fix a dimension $d\ge1$, a length $\ell>0$, and an integer $k\ge1$. On vector fields on $\T^d$, take the scalar componentwise inertia
	\begin{equation}\label{eq:inertia}
		A=(I-\ell^2\Delta)^k,\qquad u=A^{-1}m,
	\end{equation}
	and define the momentum drift
	\begin{equation}\label{eq:epdiff}
		F_A(m)=u\cdot\nabla m+(\nabla u)^{\mathsf T}m+(\nabla\cdot u)m.
	\end{equation}
	Equation~\eqref{eq:epdiff} is the EP coadjoint equation associated with the
	right-invariant metric induced by $A$ \cite{HMR}. Geometric well-posedness for
	pseudo-differential inertia is developed in \cite{BauerEtAl20}, and stochastic
	variational and analytic formulations appear in \cite{HT16,TangEP23}. The
	analysis below treats the constant-coefficient smoothing inertia
	\eqref{eq:inertia}.
	Choose
	\begin{equation}\label{eq:ep-space}
		V=H^s(\T^d;\R^d),\quad X=H^r(\T^d;\R^d),\quad
		s>d/2+2,\quad d/2+1<r<s-1,
	\end{equation}
	and
	\begin{equation}\label{eq:ep-rho}
		\rho_A(m)=\norm{m}_{W^{1,\infty}}+\norm{A^{-1}m}_{W^{1,\infty}}.
	\end{equation}
	The stochastic equation is written in the momentum variable $m$, with
	$G\in\HS(\mathfrak U,V)$. In velocity variables the corresponding additive
	coefficient is $A^{-1}G$.
	
	We then treat momentum equations with a smoothing inertia operator.
	
	\begin{proposition}[Smoothing-inertia EPDiff]\label{prop:epdiff}
		For \eqref{eq:inertia}--\eqref{eq:ep-rho}, the drift and regularizers satisfy
		Assumption~\ref{ass:local-structure} and the regularization hypotheses
		(R1)--(R4), so Theorem~\ref{thm:local-criterion} applies. Under
		Assumption~\ref{ass:stabilization},
		Theorems~\ref{thm:global},~\ref{thm:invariant}, and~\ref{thm:scaling}
		apply with the corresponding energy constant. The conserved positive quadratic
		low energy is
		\begin{equation}\label{eq:ep-energy}
			E_A(m)=\ip{m}{A^{-1}m}_{L^2}
			=\norm{A^{-1/2}m}_{L^2}^2.
		\end{equation}
		For square-root feedback and $G\ne0$, every invariant law at fixed $\eta>0$
		has infinite $H^s$ second moment.
	\end{proposition}
	\begin{proof}
		\emph{Energy form and stability in the lower topology.}
		The inverse inertia has Fourier symbol
		$\widehat{A^{-1}m}(\xi)=(1+\ell^2|\xi|^2)^{-k}\widehat m(\xi)$,
		$\xi\in\mathbb Z^d$. Hence, for every real $j$,
		\[
		\|A^{-1}m\|_{H^{j+2k}}\le C_{\ell,k}\|m\|_{H^j},\qquad
		\|u\|_{H^j}+\|\nabla u\|_{H^j}\le C_{\ell,k}\|m\|_{H^j}.
		\]
		Write $S_A(m)=(\nabla u)^{\mathsf T}m+(\nabla\cdot u)m$. The product estimate
		implies
		\begin{equation}\label{eq:ep-deformation-bound}
			\|S_A(m)\|_{H^s}
			\le C\bigl(\|\nabla u\|_\infty\|m\|_{H^s}
			+\|\nabla u\|_{H^s}\|m\|_\infty\bigr)
			\le C\rho_A(m)\|m\|_{H^s}.
		\end{equation}
		With $m_s=\Lambda^sm$, the smooth energy pairing extends as
		\begin{align}\label{eq:ep-explicit-energy-form}
			b_V(m)={}&\int_{\T^d}(\nabla\cdot u)|m_s|^2\dd x
			-2\ip{C_s(u,m)}{m_s}_{L^2}
			-2\ip{\Lambda^sS_A(m)}{m_s}_{L^2}.
		\end{align}
		It is continuous, cubic, and bounded by $C\rho_A(m)\|m\|_{H^s}^2$.
		Sobolev embedding gives the Lipschitz estimate for $\rho_A$ in $H^r$.
		
		Define the bilinear coadjoint expression
		\[
		\mathcal B_A(f,g)=(A^{-1}f)\cdot\nabla g
		+(\nabla A^{-1}f)^{\mathsf T}g+(\nabla\cdot A^{-1}f)g.
		\]
		For $w=m_1-m_2$,
		$F_A(m_1)-F_A(m_2)=\mathcal B_A(m_1,w)+\mathcal B_A(w,m_2)$.
		The transport estimate controls the energy pairing of
		$(A^{-1}m_1)\cdot\nabla w$ by $C_M\|w\|_{H^r}^2$ on $K_M$.
		For the remaining terms,
		\begin{align*}
			\|\mathcal B_A(m_1,w)-(A^{-1}m_1)\cdot\nabla w\|_{H^r}
			&\le C_M\|w\|_{H^r},\\
			\|\mathcal B_A(w,m_2)\|_{H^r}
			&\le C\|w\|_{H^r}\|m_2\|_{H^{r+1}}
			\le C_M\|w\|_{H^r}.
		\end{align*}
		These inequalities prove \eqref{eq:diff}.
		
		\emph{Regularization.}
		The same product and multiplier bounds give
		\[
		\|\mathcal B_A(f,g)\|_{H^r}\le C\|f\|_{H^r}\|g\|_{H^{r+1}},
		\qquad \|F_A(m)\|_{H^{s-1}}\le C\|m\|_{H^s}^2.
		\]
		Consequently (R2)--(R3) follow from Lemma~\ref{lem:bilinear-residual}, with
		$\delta=s-r-1$. Commutation of $A^{-1}$ and $J_n$ gives
		$\rho_A(J_nm)\le C\rho_A(m)$; self-adjointness then transfers
		\eqref{eq:ep-explicit-energy-form} to $F_n(m)=J_nF_A(J_nm)$, proving (R1).
		For (R4), Lemma~\ref{lem:transport-recovery} recovers the transport part, and
		\eqref{eq:ep-deformation-bound} gives
		\[
		\ip{J_n\Lambda^sS_A(m)}{J_nm_s}_{L^2}
		\longrightarrow\ip{\Lambda^sS_A(m)}{m_s}_{L^2}.
		\]
		The same estimates give the uniform bound in (R4), completing the verification
		of (R1)--(R4).
		
		\emph{Conserved low energy.}
		Since $u=A^{-1}m$, integration by parts gives
		\begin{align*}
			\ip{u}{u\cdot\nabla m+(\nabla\cdot u)m}_{L^2}
			&=-\ip{m}{u\cdot\nabla u}_{L^2},\\
			\ip{u}{(\nabla u)^{\mathsf T}m}_{L^2}
			&=\ip{m}{u\cdot\nabla u}_{L^2}.
		\end{align*}
		Hence $DE_A(m)F_A(m)=0$. The injectivity of $A^{-1/2}$ makes $q_0>0$ when
		$G\ne0$, and $E_A(m)\le C\|m\|_{H^s}^2$. For square-root feedback,
		\eqref{eq:upper-comparison} holds, so Theorem~\ref{thm:low} applies.
	\end{proof}
	
	\begin{remark}
		The proof uses the constant-coefficient inertia \eqref{eq:inertia} and its
		commutation with the Fourier regularizers. Variable-coefficient inertia operators
		require corresponding commutator and domain estimates.
	\end{remark}
	
	\subsection{Camassa--Holm and Degasperis--Procesi equations}\label{sec:chdp}
	On $\T$, put $A_1=1-\partial_x^2$ and
	\begin{equation}\label{eq:b-family}
		F_b(u)=u u_x+\partial_xA_1^{-1}
		\left(\frac b2u^2+\frac{3-b}{2}u_x^2\right).
	\end{equation}
	Equivalently, $m=A_1u$ obeys $m_t+u m_x+b u_xm=0$.
	The values $b=2$ and $b=3$ give the Camassa--Holm (CH) and
	Degasperis--Procesi (DP) equations. CH is the EP equation for a right-invariant
	$H^1$ metric and is bi-Hamiltonian \cite{CH,HMR,KM03}; DP is an integrable
	bi-Hamiltonian peakon equation with a distinct quadratic conserved form
	\cite{DHH,LL}. These structures place both endpoint models within the Hamiltonian framework,
	while the analytic estimates below hold for every fixed $b\in\R$.
	
	Stochastic CH equations with additive, multiplicative, transport, and jump noise
	are studied in \cite{CGG12,CDG23,CDG21global,CG16,CG17,CDG21jump,CL22,
		CDGG24,CDG24,RT,MRT,TY23,APT}. For DP, see \cite{CGDP15,ACCDP23,CCDP24,
		CGLDP25}. Recent work on long-time CH dynamics includes \cite{APT}.
	Choose
	\begin{equation}\label{eq:ch-space}
		s>\tfrac52,\quad\tfrac32<r<s-1,\quad
		V=H^s(\T),\quad X=H^r(\T),\quad\rho(u)=\norm{u}_{W^{1,\infty}}.
	\end{equation}
	
	The following proposition specializes the one-dimensional momentum framework to the
	Camassa–Holm and Degasperis–Procesi equations.
	
	\begin{proposition}[Camassa--Holm and Degasperis--Procesi]\label{prop:chdp}
		For each fixed $b\in\R$, \eqref{eq:b-family} satisfies
		Assumption~\ref{ass:local-structure} and the regularization hypotheses
		(R1)--(R4) on the spaces \eqref{eq:ch-space}. Hence
		Theorem~\ref{thm:local-criterion} gives a unique maximal $V$-valued,
		$X$-stable local solution. Under Assumption~\ref{ass:stabilization},
		Theorems~\ref{thm:global},~\ref{thm:invariant}, and~\ref{thm:scaling}
		apply with $C_V=C_{s,b}$. For $b=2$ and $b=3$,
		Theorem~\ref{thm:low} also applies,
		using the distinct conserved forms below.
	\end{proposition}
	\begin{proof}
		Write $N_b(u)=F_b(u)-u u_x$ for the nonlocal term. The order-minus-one
		multiplier and product estimates give
		\[
		\norm{N_b(u)}_{H^s}
		\le C_b\left(\norm{u}_\infty\norm{u}_{H^{s-1}}
		+\norm{u_x}_\infty\norm{u_x}_{H^{s-1}}\right)
		\le C_{s,b}\rho(u)\norm{u}_{H^s}.
		\]
		The transport-energy estimate therefore yields \eqref{eq:b}, with the factor
		$2$ absorbed in $C_{s,b}$. For differences, the transport terms are handled as
		above. The derivative in the nonlocal difference appears through
		$u_x^2-v_x^2=(u_x+v_x)(u-v)_x$.
		Since $r-1>1/2$, $H^{r-1}(\T)$ is an algebra, so the order-minus-one multiplier gives
		\[
		\norm{\partial_xA_1^{-1}[(u_x+v_x)(u-v)_x]}_{H^r}
		\le C_R\norm{u-v}_{H^r}.
		\]
		This proves \eqref{eq:diff}; the growth and feedback bounds follow likewise.
		For the approximation, polarize the drift as
		\[
		\mathcal B_b(f,g)=f g_x+\partial_xA_1^{-1}
		\left(\frac b2fg+\frac{3-b}{2}f_xg_x\right).
		\]
		Since $H^{r-1}(\T)$ is an algebra,
		$\|\mathcal B_b(f,g)\|_{H^r}\le C_b\|f\|_{H^r}\|g\|_{H^{r+1}}$,
		and $\|F_b(u)\|_{H^{s-1}}\le C_b\|u\|_{H^s}^2$.
		Lemma~\ref{lem:bilinear-residual} gives (R2)--(R3) for
		$J_nF_b(J_nu)$. Its smooth-state energy bound gives (R1). For (R4), apply
		Lemma~\ref{lem:transport-recovery} to $uu_x$; the nonlocal term $N_b(u)\in H^s$
		passes directly in the energy pairing. These formulas also define the continuous cubic extension $b_V$.
	\end{proof}
	
	For CH the conserved form is $E_{\mathrm{CH}}(u)=\norm{u}_{H^1}^2$.
	For a general $b$, direct integration by parts gives
	\begin{equation}\label{eq:b-cancel}
		\ip{A_1u}{F_b(u)}_{L^2}
		=\int u\,[u m_x+b u_xm]\dd x
		=\frac{b-2}{2}\int u_x^3\dd x.
	\end{equation}
	Thus the $H^1$ cancellation is specific to $b=2$.
	
	For DP take $B_4=4-\partial_x^2$ and
	\begin{equation}\label{eq:dp-energy}
		K_{\mathrm{DP}}=A_1B_4^{-1},\qquad
		E_{\mathrm{DP}}(u)=\ip{u}{K_{\mathrm{DP}}u}_{L^2}
		=\sum_{n\in\mathbb Z}\frac{1+n^2}{4+n^2}|\widehat u(n)|^2.
	\end{equation}
	The multiplier lies in $[1/4,1]$, so $E_{\mathrm{DP}}$ is equivalent to the
	$L^2$ norm squared. This positive quadratic energy is standard in DP analysis
	\cite{LL}. Its periodic cancellation follows directly: set $v=B_4^{-1}u$, so
	$K_{\mathrm{DP}}u=u-3v$ and
	$F_3(u)=u u_x+\tfrac32A_1^{-1}\partial_x(u^2)$. Then
	\begin{align}\label{eq:dp-cancel}
		\ip{K_{\mathrm{DP}}u}{F_3(u)}
		&=\int(u-3v)u u_x\dd x+\frac32\int v\partial_x(u^2)\dd x\notag\\
		&=\frac32\int v_xu^2\dd x-\frac32\int v_xu^2\dd x=0.
	\end{align}
	For nonzero $G$, both the CH and DP quadratic energy injection rates are positive. Their low energies are bounded by $C\norm{u}_{H^s}^2$, giving infinite stationary $H^s$ second moments for square-root feedback.
	
	\subsection{Hunter--Saxton equation}\label{sec:hs}
	The Hunter--Saxton equation is the EP equation for the homogeneous
	$\dot H^1$ metric on the diffeomorphism group modulo rotations and admits a
	bi-Hamiltonian formulation \cite{HS,KM03}. Its stochastic characteristic theory
	is developed in \cite{HKP21}. We use a periodic mean-zero momentum formulation
	with radial It\^o noise. In one dimension set $m=-u_{xx}$ and let
	$(-\partial_x^2)^{-1}$ act on the mean-zero subspace. The momentum equation is
	\[
	m_t+u m_x+2u_xm=0.
	\]
	
	We finish the model verification with the mean-zero Hunter–Saxton equation.
	
	\begin{proposition}[Mean-zero Hunter--Saxton equation]\label{prop:hs}
		Let
		\[
		V=H_0^s(\T),\qquad X=H_0^r(\T),\qquad
		s>\tfrac52,\quad \tfrac32<r<s-1,
		\]
		and assume that the noise takes values in $V$. With
		$\rho(m)=\|m\|_{W^{1,\infty}}+\|(-\partial_x^2)^{-1}m\|_{W^{1,\infty}}$,
		the drift satisfies Assumption~\ref{ass:local-structure} and the regularization
		hypotheses (R1)--(R4). Under Assumption~\ref{ass:stabilization},
		Theorems~\ref{thm:global},~\ref{thm:invariant}, and~\ref{thm:scaling}
		apply. The conserved positive low energy is
		\[
		E_{\mathrm{HS}}(m)
		=\ip{m}{(-\partial_x^2)^{-1}m}_{L^2}
		=\|u_x\|_{L^2}^2.
		\]
	\end{proposition}
	\begin{proof}
		The inverse gains two derivatives on the mean-zero subspace, so the estimates in
		Proposition~\ref{prop:epdiff} apply with $A=-\partial_x^2$.  The drift preserves
		mean zero because
		\[
		\int_{\T}(u m_x+2u_xm)\dd x
		=\int_{\T}u_xm\dd x
		=-\int_{\T}u_xu_{xx}\dd x=0.
		\]
		The Fourier regularizers preserve the constraint and commute with the inverse.
		More explicitly, $A^{-1}$ has symbol $n^{-2}$ for $n\ne0$ and zero symbol at
		$n=0$, so $\|A^{-1}m\|_{H^{j+2}}\le C\|m\|_{H^j}$ on the constrained
		space. The residual decomposition \eqref{eq:three-residuals} applies to
		$\mathcal B_A(f,g)=(A^{-1}f)g_x+2(\partial_xA^{-1}f)g$ in the ambient
		Sobolev spaces. The full diagonal drift preserves mean zero by the displayed
		integral identity. Self-adjointness gives (R1),
		Lemma~\ref{lem:bilinear-residual} gives (R2)--(R3), and
		Lemma~\ref{lem:transport-recovery} gives (R4) for the transport term;
		$2u_xm\in H^s$ is a regular remainder. The low-energy identity is the EP
		cancellation on the quotient by constants. Both the drift and the noise preserve
		the mean-zero constraint.
	\end{proof}

	\medskip 
	
	\appendix
	
	\section{A fixed-forcing moment obstruction}\label{sec:low}
	At fixed $\eta>0$, an additional conserved positive quadratic energy yields an
	obstruction to finite high-order second moments. The quadratic form may be a
	Hamiltonian, a Casimir, or another conserved energy.
	
	Let $\mathcal{H}_0$ be a Hilbert energy space and $T:V\to\mathcal{H}_0$ a continuous injective linear map. Write
	\begin{equation}\label{eq:lowenergy}
		E_0(u)=\norm{Tu}_{\mathcal{H}_0}^2,\quad
		q_0=\norm{TG}_{\HS(\mathfrak U,\mathcal{H}_0)}^2,\quad
		S_0(u)=\norm{(TG)^*Tu}_{\mathfrak U}^2\le q_0E_0(u).
	\end{equation}
	This formulation includes a positive quadratic form $\ip{Ku}{u}$ by taking
	$T=K^{1/2}$ on its domain and accommodates unbounded forms through their natural
	domains.
	
	\begin{assumption}[Conserved low energy]\label{ass:low}
		Along the local $V/X$-solutions, the low-energy identity holds with no
		deterministic drift contribution:
		\begin{align}
			\dd E_0(U_\eta)&=[a(U_\eta)^2E_0(U_\eta)+\eta^2q_0]\dd t+\dd M_t,
			\label{eq:lowito}\\
			\dd[M]_t&=4[a(U_\eta)^2E_0(U_\eta)^2+\eta^2S_0(U_\eta)]\dd t.
			\label{eq:lowqv}
		\end{align}
		For smooth states this requires $DE_0(u)F(u)=0$. We assume $q_0>0$.
	\end{assumption}
	
	For the conserved energies verified in Section~\ref{sec:models}, the map $T$
	extends to a bounded linear map $X\to\mathcal H_0$, where $\mathcal H_0$ is
	the scalar or vector $L^2$ space, or the product $L^2$ space for MHD\@. The
	corresponding choices are
	\[
	T=\begin{cases}
		I, & \text{active scalars and ideal MHD},\\
		A^{-1/2}, & \text{EPDiff},\\
		A_1^{1/2}, & \text{Camassa--Holm},\\
		K_{\mathrm{DP}}^{1/2}, & \text{Degasperis--Procesi},\\
		(-\partial_x^2)^{-1/2}, & \text{mean-zero Hunter--Saxton}.
	\end{cases}
	\]
	The last operator acts only on the mean-zero subspace. For CH, the required
	boundedness follows from $r>3/2>1$; the other displayed multipliers have
	nonpositive order. Each map is injective on its solution space $V$.
	The smooth-state cancellations extend by density and continuity of
	$F:V\to X$, giving
	\begin{equation}\label{eq:low-cancellation-extension}
		2\ip{Tu}{TF(u)}_{\mathcal H_0}=0,\qquad u\in V.
	\end{equation}
	To justify the stochastic identity, stop the $V/X$-solution at
	$\tau_R=\inf\{t\ge0:\|U_\eta(t)\|_V\ge R\}$, with an additional localization on the initial norm when the random initial
	state is unbounded. Applying the bounded map $T$
	to the equation in $X$ yields an $\mathcal H_0$-valued semimartingale:
	\begin{align}\label{eq:low-transformed-equation}
		TU_\eta(t\wedge\tau_R)
		={}&TU_\eta(0)-\int_0^{t\wedge\tau_R}TF(U_\eta(s))\dd s\notag\\
		&+\int_0^{t\wedge\tau_R}a(U_\eta(s))TU_\eta(s)\dd W_s
		+\eta\int_0^{t\wedge\tau_R}TG\dd B_s.
	\end{align}
	All coefficients in this stopped equation are locally integrable in the
	corresponding Hilbert norms. The ordinary Hilbert-space It\^o formula for
	$\|TU_\eta\|_{\mathcal H_0}^2$ is therefore applicable. By
	\eqref{eq:low-cancellation-extension}, its finite-variation part is
	$[a(U_\eta)^2E_0(U_\eta)+\eta^2q_0]\dd t$, and its local martingale part is
	\[
	M_t=2\int_0^t a(U_\eta(s))E_0(U_\eta(s))\dd W_s
	+2\eta\int_0^t
	\ip{(TG)^*TU_\eta(s)}{\dd B_s}_{\mathfrak U},
	\]
	with the same stopping understood. Independence of $W$ and $B$ gives
	\eqref{eq:lowqv}. Thus \eqref{eq:lowito}--\eqref{eq:lowqv} hold along the stochastic solutions.
	Finally, $G\ne0$ and injectivity of $T$ imply $q_0>0$, completing the
	model-specific verification of Assumption~\ref{ass:low}.
	
	The first appendix result describes a limitation on higher stationary moments at fixed
	positive forcing.
	
	\begin{theorem}[A stationary weighted-moment obstruction]\label{thm:low}
		Under the hypotheses of Theorem~\ref{thm:invariant} and
		Assumption~\ref{ass:low}, every $\mu\in\cI_\eta$, for each fixed
		$\eta\in(0,1]$, satisfies
		\begin{equation}\label{eq:weighted-divergence}
			\int_V a(u)^2\sqrt{1+E_0(u)}\,\mu(\dd u)=\infty.
		\end{equation}
		If in addition
		\begin{equation}\label{eq:upper-comparison}
			a(u)^2\sqrt{1+E_0(u)}\le C(1+\norm{u}_V^2),
		\end{equation}
		then $\int\norm{u}_V^2\mu(\dd u)=\infty$. Condition~\eqref{eq:upper-comparison} holds for square-root feedback and a
		continuous quadratic energy controlled by the $V$ norm. For faster-growing
		coefficients $\psi$, it must be verified separately.
	\end{theorem}
	\begin{proof}
		Let $U_\eta$ be a stationary solution with one-time law $\mu$. Put
		$Y=E_0(U_\eta)$ and $R_0=\sqrt{1+Y}$. Applying It\^o's formula to the stopped
		process gives
		\begin{equation}\label{eq:root-drift}
			\dd R_0=D(U_\eta)\dd t+\dd N_t,\qquad
			D(u)=\frac{a(u)^2Y+\eta^2q_0(1+Y)-\eta^2S_0(u)}{2(1+Y)^{3/2}}.
		\end{equation}
		Because $S_0\le q_0Y$ and $\eta^2q_0>0$, $D(u)>0$ for every $u$. Moreover,
		\begin{equation}\label{eq:root-bracket}
			\frac{\dd[N]_t}{\dd t}
			=A(U_\eta):=\frac{a(U_\eta)^2Y^2+\eta^2S_0(U_\eta)}{1+Y}
			\le a(U_\eta)^2R_0^2+\eta^2q_0.
		\end{equation}
		Suppose, for contradiction, that $\int a^2R_0\dd\mu<\infty$.
		Then $a^2$ is integrable and, since $a^2\ge b_0$, so is $R_0$.
		The drift $D$ is integrable because
		$0<D\le(a^2+\eta^2q_0)/(2R_0)$.
		
		Choose bounded increasing concave functions $\Theta_n$ with
		$\Theta_n'=1$ below $n$, $\Theta_n'=0$ above $2n$, and
		$|\Theta_n''|\le C/n$ supported on $[n,2n]$.
		For fixed $n$, the martingale in the formula for $\Theta_n(R_0)$ is square integrable: on $R_0\le2n$, its bracket is bounded by
		$4n^2a^2+\eta^2q_0$, which is integrable under the contradiction hypothesis. Localization can therefore be removed. Stationarity yields
		\begin{equation}\label{eq:root-stationary}
			0=\int\Theta_n'(R_0)D\dd\mu
			+\frac12\int\Theta_n''(R_0)A\dd\mu.
		\end{equation}
		The absolute value of the second term is at most
		\[
		C\int_{\{R_0\ge n\}}a^2R_0\dd\mu+C\eta^2q_0/n\longrightarrow0.
		\]
		Dominated convergence in the first term gives $\int D\dd\mu>0$, a
		contradiction. This proves \eqref{eq:weighted-divergence}; the final assertion
		follows from \eqref{eq:upper-comparison}.
	\end{proof}
	
	\begin{remark}
		The condition $q_0>0$ requires energy injection in at least one direction. The
		conclusion concerns fixed nonzero forcing and establishes moment divergence; the
		small-forcing Cauchy-law limit describes a separate asymptotic regime.
	\end{remark}

	\medskip 
	
	\section{Linear limits with a retained deterministic mode}\label{sec:linear}
	The main scaling theorem concerns a drift that vanishes after rescaling. This
	appendix treats the linear equation obtained when a deterministic mode remains.
	Its covariance formula also clarifies the purely radial limit in
	Theorem~\ref{thm:linear-limit}.
	
	Let $-\mathscr A$ generate a strongly continuous semigroup $S(t)=e^{-t\mathscr A}$ on a real separable Hilbert space $X$ and suppose
	\begin{equation}\label{eq:semigroup}
		\norm{S(t)}_{\mathcal L(X)}\le M_Se^{\chi t},\qquad t\ge0.
	\end{equation}
	Consider the mild equation
	\begin{equation}\label{eq:linear}
		\dd V+\mathscr A V\dd t=a_0V\dd W_t+G\dd B_t,
		\qquad G\in\HS(\mathfrak U,X),
	\end{equation}
	and set $Q=GG^*$, where the adjoint is taken in $X$.
	
	The final result records the stationary law when a finite-dimensional linear Hamiltonian
	mode is retained.
	
	\begin{theorem}[Stationary law with a retained linear mode]\label{thm:linear}
		If $\chi<a_0^2/2$, \eqref{eq:linear} has a unique invariant probability law
		$\nu_{\mathscr A}$. This law is a mixture of centered Gaussian measures with
		random covariance
		\begin{equation}\label{eq:random-cov}
			Q_W=\int_0^\infty e^{2a_0\widetilde W_r-a_0^2r}
			S(r)QS(r)^*\dd r.
		\end{equation}
		The integral is almost surely trace class. Set
		$\lambda_\alpha=\alpha(1-\alpha)a_0^2/2$. For
		\begin{equation}\label{eq:lin-q}
			0<\alpha<1,\qquad
			\lambda_\alpha-\alpha\chi>0,
		\end{equation}
		the law has a finite $\alpha$th norm moment, and the common-noise coupling
		satisfies
		\begin{equation}\label{eq:lin-contraction}
			\E\norm{V_t^x-V_t^y}_X^\alpha
			\le M_S^\alpha e^{-(\lambda_\alpha-\alpha\chi)t}\norm{x-y}_X^\alpha.
		\end{equation}
		If
		\begin{equation}\label{eq:cov-invariance}
			S(t)QS(t)^*=Q,\qquad t\ge0,
		\end{equation}
		then $\nu_{\mathscr A}=\Law(Z/(a_0|\zeta|))$, where $Z$ has covariance
		$Q$.
	\end{theorem}
	\begin{proof}
		The mild finite-time formula is
		\[
		V_t=\mathcal E_tS(t)V_0+
		\mathcal E_t\int_0^t\mathcal E_s^{-1}S(t-s)G\dd B_s.
		\]
		Finite-past solutions give a stationary conditional Wiener integral as in
		\eqref{eq:stationary-convolution}, with the factor $S(t-s)$ inserted.
		Its conditional covariance is \eqref{eq:random-cov}. The trace is bounded by
		\[
		M_S^2\norm{G}_{\HS(\mathfrak U,X)}^2
		\int_0^\infty e^{2a_0\widetilde W_r-(a_0^2-2\chi)r}\dd r,
		\]
		which is finite almost surely because $\widetilde W_r/r\to0$ and
		$a_0^2-2\chi>0$. Taking the finite-past limit and using the mild flow identity
		yields a stationary solution. Conditionally on $W$, the infinite-past expression
		is a Wiener integral with respect to $B$; equivalently, it is the limit of
		finite-past adapted solutions.
		
		To obtain the $\alpha$th moment, split the last scalar integral into unit
		intervals and use subadditivity of the power $\alpha/2<1$. After taking that
		power, the expected contribution from $[n,n+1]$ is bounded by a constant times
		\[
		\E e^{\alpha a_0\widetilde W_n-\alpha(a_0^2/2-\chi)n}
		=e^{-(\lambda_\alpha-\alpha\chi)n}.
		\]
		Brownian exponential suprema on a unit interval have finite moments. The geometric series converges under \eqref{eq:lin-q}. Conditional Gaussian moments are bounded by a constant times $(\Tr Q_W)^{\alpha/2}$, proving the claim.
		The difference formula is $V_t^x-V_t^y=\mathcal E_tS(t)(x-y)$, which proves
		\eqref{eq:lin-contraction}. This difference also tends to zero almost surely. A
		common-noise coupling and bounded convergence then give uniqueness among all
		invariant probabilities. Under \eqref{eq:cov-invariance}, $Q_W=I\,Q$, where $I$ has the distribution \eqref{eq:dufresne}; this proves the final statement.
	\end{proof}
	
	\begin{remark}[Infinite-dimensional Gaussian notation]
		The mixture is defined by
		$\nu_{\mathscr A}(A)=\E[\mathcal N(0,Q_W)(A)]$. In infinite dimensions, the
		identity operator is not trace class, so an $X$-valued Gaussian with identity
		covariance does not exist. The notation $Q_W^{1/2}\xi$ therefore requires a
		cylindrical Gaussian $\xi$ or an equivalent trace-class series construction.
	\end{remark}
	
	\begin{example}[Covariance compatibility]\label{ex:rotation}
		Take $X=\R^2$, $S(t)$ a nontrivial rotation group, and
		$Q=e_1e_1^{\mathsf T}$. For every $h\ne0$, the quadratic form of
		\eqref{eq:random-cov} is
		\[
		\int_0^\infty e^{2a_0\widetilde W_r-a_0^2r}
		|h\cdot S(r)e_1|^2\dd r>0.
		\]
		Thus $Q_W$ is positive definite almost surely and the mixture has full
		support. By contrast, a ratio law built from the original rank-one covariance
		$Q$ is supported on a line. Hence covariance compatibility is necessary for the
		explicit ratio representation.
	\end{example}

	\medskip 
	
	\noindent\textbf{Funding.}
	This work was supported by the Baima Lake Laboratory Joint Fund of Zhejiang Provincial NSF of China [Grant No. LBMHZ26E060006], the National Natural Science Foundation of China [Grant Nos. 12571164, 12671292], the Jiangsu Provincial Scientific Research Center of Applied Mathematics [Grant No. BK20233002], and the Fundamental Research Funds of Zhejiang Sci-Tech University [Grant No. 26062140-Y].
	
	\medskip
	\noindent\textbf{Data availability.}
	No datasets were generated or analyzed for this study.
	
	\medskip
	\noindent\textbf{Declaration of AI-assisted tools.}
During the preparation of this work, the authors used ChatGPT(OpenAI) to assist with language editing, manuscript organization, and checking the consistency of terminology and notation. After using this tool, the authors reviewed and edited the content as needed, and take full responsibility for the content of article.

\end{document}